\input amstex
\documentstyle{amsppt}

\hsize=4.75in
\vsize=8in

\def\J{{\Cal J}}
\def\I{{\Cal I}}

\def\O{{\Cal O}}
\def\H{{\Cal H}}
\def\Aut{\hbox{\rm Aut}}
\def\K{{\Cal K}}

\def\T{{\Cal T}}

\def\F{{\Cal F}}
\def\L{{\Cal L}}

\def\A{{\Cal A}}
\def\B{{\Cal B}}

\def\Z{{\Cal Z}}
\def\D{{\Cal D}}
\def\ov{\overline}

\rightheadtext {$C^*$--algebras generated by Hilbert bimodules}
\leftheadtext {T.~Kajiwara, C.~Pinzari, Y.~Watatani}
\topmatter
\title
Ideal structure and simplicity\\ of the $C^*$--algebras generated by 
Hilbert bimodules\\ \\ 
({\it Revised version, January 1998})\\ 
\endtitle

\author 
Tsuyoshi Kajiwara \footnote{Supported by the
Grants--in--aid for Scientific Research, The Ministry of Education,
Science and Culture, Japan} Claudia Pinzari\footnote{Supported by MURST,
CNR--GNAFA and European Community} and Yasuo Watatani$^1$
\endauthor

\affil
Department of Environmental and Mathematical Sciences, \\
         Okayama University, Tsushima, 700, Japan\\ \\
Dipartimento di Matematica, \\
         Universit\`a di Roma Tor Vergata, 00133 
         Roma, Italy\\ \\
Graduate School of Mathematics, \\
         Kyushu University, Ropponmatsu, Fukuoka, 810, Japan
\endaffil

\abstract
   Pimsner introduced the $C^*$--algebra $\Cal O_X$ generated by a Hilbert bimodule 
$X$ over a $C^*$--algebra $\A$.   
   We look for additional conditions that $X$ should satisfy in order
to study the simplicity and, more generally, the ideal
structure of   $\O_X$ when $X$ is finite projective. 
We introduce two conditions: ``$(I)$--freeness'' and ``$(II)$--freeness'', stronger than 
the former, in analogy with 
 [J. Cuntz, W. Krieger, Invent. Math. 56, 251--268] and 
[J. Cuntz, Invent. Math. 63, 25--40]
respectively.  $(I)$--freeness comprehend the case of
 the bimodules associated with an inclusion of simple $C^*$--algebras
with finite index, real or pseudoreal bimodules with finite dimension 
 and the case of ``Cuntz--Krieger bimodules''. If $X$ satisfies this condition
 the $C^*$--algebra $\O_X$ does not depend on the choice
of the generators when $\A$ is faithfully represented.
     
   As a consequence, if $X$ is $(I)$--free and $\A$ is $X$--simple, then $\O_X$ is
simple.  In the case of Cuntz--Krieger algebras $\O_A$,
$X$--simplicity corresponds to the irreducibility of the
matrix $A$. If  $\A$ is simple and p.i. then $\O_X$ is p.i., if $\A$ is nonnuclear then $\O_X$
is nonnuclear. We therefore  provide  examples of (purely) infinite 
nonnuclear simple $C^*$--algebras.

Furthermore if $X$ is $(II)$--free, we determine the
ideal structure of $\O_X$.
\endabstract

\endtopmatter

\document

\heading
1. Introduction
\endheading

  In \cite{Pi}  Pimsner  associated with 
every Hilbert $C^*$--bimodule  $X$ a $C^*$--algebra $\Cal O_X$. The bimodule
$X$ has a natural  
 representation in $\O_X$, and generates $\O_X$ as a $C^*$--algebra.
 As noted in \cite{Pi},
  representing Hilbert $C^*$--bimodules in $C^*$--algebras
is a natural generalization of the
 notion of  Hilbert space contained in a $C^*$--algebra,
introduced by Doplicher and Roberts (see, e.g., \cite{DR2}).

The construction of $\O_X$ generalizes and unifies various  examples of
$C^*$--algebras. As noted in \cite{Pi},
 $\Cal O_X$ is a generalization of  Cuntz--Krieger algebras
$\Cal O_A$ \cite{Cu, CK},  crossed product $C^*$--algebras
$\A \rtimes \Bbb Z$ by an automorphism and 
crossed products by
partial automorphisms introduced by  Exel \cite{Ex}.
       In the particular case where $X$ is the bimodule
associated with an inclusion of von Neumann algebras, the 
same construction had been  independently considered by Katayama
\cite{Ka} .
If $X$ is a Hilbert $C^*$--bimodule 
in the sense of \cite{BMS}   
  $\Cal O_X$,  previously considered in \cite{AEE} from the viewpoint
of crossed product $C^*$--algebras, was introduced as a    ``crossed product 
by $X$''. Furthermore,  in the case where
$X$ is full and finite projective as a right $\A$--module,    $\O_X$ is  the 
  Doplicher--Roberts 
algebra  associated to $X$ regarded as an object of the
semitensor $C^*$--category with objects Hilbert $C^*$--bimodules
over $\A$ and arrows right adjointable maps (cf. \cite{DPZ} or Proposition 2.5).
  
Examples of Hilbert $C^*$--bimodules are provided from 
inclusions of simple $C^*$--algebras with finite index studied in \cite{Wa} and 
Izumi \cite{I1}, \cite{I2}.

   Our main interest in this paper is to study  the ideal structure
of   $\Cal O_X$ when $\A$ is unital and $X$ is finite projective. 
We shall discuss the case where $\A$ is $\sigma$--unital
and $X_{\A}$ is countably generated elsewhere \cite{P, KWP}. 

We look for a general
criterion that generalizes and unifies the case  of crossed products by
properly outer automorphisms [El, OP1, Co1, Co2, CT, Ol, Ki, OP2], condition 
$(I)$ and $(II)$ for Cuntz--Krieger algebras
[Cu, CK, Cu2],  condition of having ``determinant 1'' for Doplicher and Roberts 
algebras [DR1, DR2, DR3],
and the case of bimodules associated 
with inclusions of simple $C^*$--algebras with
finite Jones index $>1$ [Jo, Wa].

 In \S 4 we formulate a condition,
that we call $(I)$--freeness, in analogy with   condition $(I)$ in \cite{CK}.
Our main general result is Theorem 4.3.
We prove that if $X$ is $(I)$--free then any representation of the bimodule
 $X$ in a $C^*$--algebra $\D$ that is faithful on the coefficient algebra $\A$
extends to a {\it faithful} representation of $\O_X$ in $\D$. Hence $\O_X$
is the unique $C^*$--algebra generated by a 
faithful representation of $X$ in a $C^*$--algebra.
In particular, if $\A$ is $X$--simple, $\O_X$ is simple.

 Furthermore we introduce
a stronger condition, $(II)$--freeness, again in analogy with property $(II)$ of
\cite{Cu2}, that allows us to determine closed ideals of $\O_X$ in terms of a certain
class of closed ideals of $\A$. 

The strategy of the proof of Theorem 4.3 is the following. 
First we identify
the algebra $\O_X$ with the Doplicher--Roberts algebra associated with $X$
regarded as an object of the right tensor $C^*$--category of Hilbert $C^*$--bimodules
with arrows right adjointable maps (\S 2). Thus   $\O_X$ carries a natural
structure of ${\Bbb Z}$--graded $C^*$--algebra $\oplus_{k\in{\Bbb Z}}{\O_X}^{(k)}$,
similar to that of Cuntz--Krieger algebras and crossed product $C^*$--algebras
by an automorphism.

We next  need show that the canonical dense $^*$--subalgebra $^0{\O_X}$ generated by
the homogeneous spaces has a unique $C^*$--norm and that any representation
of $\O_X$ that is faithful on $X$ is faithful. To this aim it is enough
to show, by a result of \cite{DR1},
 that the projection of $^0{\O_X}$ onto ${\O_X}^{(0)}$ is continuos in
any $C^*$--norm. This second step is present, for example, in \cite{El, Ki} for
certain crossed product $C^*$--algebras by outer actions of
automorphism groups, in \cite{Cu, CK} for
the algebras $\O_n$ and $\O_A$, where $A$ is a matrix 
 that satisfies the so called condition $(I)$, and in \cite{DR1}  for
  the algebras $\O_\rho$
associated with a unitary representation of a compact group with determinant $1$.
We introduce  our main assumption that $X$ is $(I)$--free in order to prove this 
second step.

It is nice to note that if $X$ is $(I)$--free,
the simplicity criterion that we get, that requires for $\A$ to be $X$--simple,
(i.e. lack of
the so called $X$--invariant ideals of $\A$), corresponds, in the case of Cuntz--Krieger
algebras $\O_A$ to irreducibility of the matrix $A$.

Thirdly, in order to study the ideal structure of $\O_X$ we use the condition 
that $X$ is $(II)$--free. This condition, similarly to condition $(II)$ by
Cuntz and Krieger \cite{CK2}, assures that any closed ideal $\J$  of
$\O_X$ is ${\Bbb Z}$--graded, and hence determined by  a certain ideal
of $\A$, and it  guarantees that $\O_X/\J$ is still
 an algebra of the form
$\O_Y$, where $Y$ is a $(I)$--free Hilbert $C^*$--bimodule obtained as a quotient of
$X$ over a quotient of $\A$.

In \S 5 we discuss two applications of our main theorem: the case of 
real or pseudoreal bimodules with Jones index $>1$, that contains the special
important case of bimodules associated with proper inclusions of simple 
$C^*$--algebras with finite index in the sense of
\cite{Wa}, and the case of ``Cuntz--Krieger'' bimodules, 
 namely
bimodules for which the coefficient
algebra
is the direct sum of $d>1$ simple $C^*$--algebras and furthermore the 
adjacency matrix satisfies condition $(I)$ of \cite{CK}.

  In \S 6 we investigate about pure infiniteness of 
$\O_X$.
We prove that if $X$ is $(I)$--free and $\A$ is simple and purely infinite
then $\O_X$ is simple  and purely infinite.
   There exists  a relation with a classification of purely infinite
simple $C^*$--algebras.  We do not pursue this in the present paper,
just see, for example,  \cite{R\o1}, \cite{R\o2} and other works of recent
rapid development.

After this work was completed we received the paper
\cite{MS} by  Muhly and  Solel,
where similar problems  have been considered.

\heading
2. Representations of  Bimodule Structures
\endheading

	Following  \cite{Pi},  a Hilbert bimodule $X$ over a $C^*$--algebra $\A$ is
 a right Hilbert $\A$--module (with $\A$--valued inner product denoted by
$(\ \vert \ )_\A$)
   endowed with an
isometric $^*$--homomorphism $\phi : \A \rightarrow \Cal L_\A(X)$ from $\A$ to the algebra
of  right $\A$--linear adjointable maps on $X$ (the left action of $\A$).
 In the following, we will write  $X_\A$ when
 $X$ is regarded  just as a  Hilbert $\A$--module,  and $_\A X_\A$
when it is regarded as a Hilbert $\A$--bimodule.

    In this section we review the various constructions  of the
$C^*$--algebra $\Cal O_X$ following  
 \cite{Pi}, \cite{AEE}, \cite{DPZ} and \cite{Ka}. We will assume here and in the rest of 
this paper,
 for simplicity, that the $C^*$--algebra $\A$  is
unital and that the Hilbert $C^*$--module $X_\A$ is  full and finite projective.

      For $x,y \in X$, the rank one operator $\theta _{x,y} \in \Cal L_\A(X_\A)$
is defined by $\theta _{x,y}(z) = x(y\vert z)_\A$ for $z \in X$.
   We denote by $\Cal K_\A(X_\A)$ the norm closure of linear combinations of
rank one operators.    A finite subset $\{u_1,\dots,u_n\} \subset X$ is called a 
{\it basis} for $X$
if $x = \sum_{i=1}^n u_i(u_i\vert x)_\A$  \ for all $x \in X$.
Recall that since $X_\A$ is finite projective, it has a basis, and that
 every adjointable map on $X_\A$ is an element of $\K_\A(X_\A)$. 
We set $\K:=\K_\A(X_\A)=\L_\A(X_\A)$.
       
   We fix a  basis $\{u_1,\dots,u_n\}$ of $X_\A$.
Let $\phi : \A \rightarrow \Cal L_\A(X_\A)$ be the defining  left action 
of $\A$ on $X$.
   For $a \in \A$, 
$$
  \phi(a)u_j = \sum_{i=1}^n u_ia_{ij}
$$
where $a_{ij} := (u_i\vert \phi(a)u_j)_\A \in \A$.
   Then the $C^*$--algebra $ O_X$ is the universal $C^*$--algebra generated
by $\A$ and n operators $\{S_1,\dots,S_n\}$ satisfying  relations
$$
   S_i^*S_j = (u_i\vert u_j)_\A, \quad
   \sum_{i=1}^n S_iS_i^* = I \quad \text{and} \quad
   aS_j = \sum_{i=1}^n S_ia_{ij} \quad \quad \text{for} \quad a \in \A,
$$
for $j = 1,\dots,n$.
   The generators $S_i (i= 1, \ldots ,n)$ are partial isometries
if and only if the basis satisfies  $(u_i \vert u_j)_\A = 0 \ (i \ne j)$,
because the orthogonality of the basis implies that 
$u_i = u_i(u_i\vert u_i)_\A$.  

 The following formulation of the  universality of  the algebra $\O_X$, due to
Pimsner
 \cite{Pi}, makes  use only of the bimodule structure of $X$, with no need of
choosing a basis of $X$.

   A representation  of  the Hilbert  $\A$--bimodule $_\A X_\A$ in a $C^*$--algebra
$\D$ is a pair $(V, \rho_\A)$ constituted by 
linear map $V : X \rightarrow \D$ and
a unital $^*$--homomorphism  
$\rho_\A : \A \rightarrow \D$ satisfying the following relations:
$$
V_{xa} = V_x\rho_\A(a),\eqno(2.1)
$$
$$
V_{\phi(a)x} = \rho_\A(a)V_x,\eqno(2.2)
$$
$$
V_x^*V_y = \rho_\A((x\vert y)_\A)\eqno(2.3)
$$
for $x, y \in X,  a \in \A$. Note that $\rho_\A$ is 
uniquely determined by $V$ since $X$ is assumed to be full. 
If $\{u_i, i=1,\dots, n\}$ is a basis
of $X$,
$\sum_{i=1}^{n} V_{u_i}{V_{u_i}}^*$ does not depend on the choice of the basis,
and is called the {\it support\/} of
$(V, \rho_\A)$. 

   We regard  $X$ 
as an imprimitivity bimodule ${}_\K X_{\A}$ in the sense of \cite{Ri},
with  left inner product ${}_\K(x\vert y) = \theta_{x,y}$ for $x,y \in X$.
   A representation  of $_\K X_\A$ in a $C^*$--algebra
$\D$ is a triple $(V, \rho_\K, \rho_\A)$ where $(V, \rho_\A)$ is a representation 
of $_\A X_\A$ in $\D$ and  $\rho_\K : \K \rightarrow \D$
is a unital $^*$--homomorphism   such that
$$
V_{kx} = \rho_\K(k)V_x, \eqno(2.4)
$$
$$
V_xV_y^* = \rho_\K({}_\K(x\vert y))\eqno(2.5)
$$
for $x, y\in X,  k \in \K$.
 Note that   since $\rho_\K$ is unital, 
$$\sum_{i=1}^{n} V_{u_i}{V_{u_i}}^*=I\ ,\eqno(2.6)$$
for any basis $\{u_i, i=1,\dots, n\}$ of $X$. Furthermore,
comparing $(2.2)$ and $(2.4)$, $\rho_\K\circ\phi=\rho_\A$.

  The $C^*$--algebra $\Cal O_X$ is the universal $C^*$--algebra
generated by operators 
$S_x$, $x \in X$ where $S: X\to\O_X$ defines a representation of the imprimitivity
bimodule  $_\K X_\A$ in 
$\O_X$.

   The equivalence of the above two definitions depends on the fact that 
any representation of the Hilbert $\A$--bimodule $_\A X_\A$ in a unital
$C^*$--algebra $\D$ with support $I$ extends to a representation of
the imprimitivity bimodule $_\K X_\A$ \cite{Pi}.
    We will give  a proof in  Lemma  2.2 and  Corollary 2.3,
because it is one of the  steps of the proof of the simplicity
of $\Cal O_X$.

\cite{Pi} and \cite{Ka}, generalizing to Hilbert $C^*$--bimodules a well known
construction for Hilbert spaces, provide a 
 concrete realization of $\Cal O_X$ as a quotient
of the Toeplitz algebra $\Cal T_X$ acting on the Fock space of $X$. This construction
shows in particular that every Hilbert $C^*$--bimodule can be {\it faithfully\/}
represented in some $C^*$--algebra.
 
   We consider right Hilbert $\A$--modules $X^{\otimes m} = X \otimes _\A \dots \otimes
_\A X$ (m-times), with the convention $X^{\otimes 0} = \A$ 
   (since $X_\A$ is finite projective, the algebraic tensor products over $\A$ are 
right Hilbert $\A$--modules, with no  need of completion (cf. \cite{KW1}).
   Let $F(X) = \oplus_{m=0}^\infty X^{\otimes m}$ be  the full Fock space of $X$,
regarded as a
 right Hilbert $\A$--module. We consider the representation of $_\A X_\A$ in 
$L_\A(F(X)_\A)$ by letting $\A$ act in the obvious way on $F(X)$ and defining,
   for $x \in X$,     the creation
operator $T_x \in \Cal L_\A(F(X)_\A)$
$$
\align
& \quad \quad T_x(x_1\otimes \dots\otimes x_m) = x\otimes x_1\otimes\dots\otimes x_m\ ,
\\
& \quad \quad T_x(a) = xa 
\endalign
$$
for $x_1\otimes\dots\otimes x_m \in X^{\otimes m}$  and $a \in \A$.
  The $C^*$--subalgebra $\Cal T _X$ of $\Cal L_\A(F(X)_\A)$
 generated by $\{T_x ; x \in X\}$ is
called the Toeplitz algebra of $X$. 

Let $\pi  : \Cal L_\A(F(X)_\A) \rightarrow
      \Cal L_\A(F(X)_\A)/\Cal K_\A(F(X)_\A)$ be the quotient map,
and set $S_x=\pi(T_x), x\in X$. Let $C^*(S_X)$ denote the $C^*$--algebra 
generated by $\{S_x ; x \in X\}$.   For $x = (x_1,\dots,x_r) \in X^{\times r}$, we set
$S_x = S_{x_1}\dots S_{x_k}$ and
  put $S_{\emptyset} = I$.
  Let $\Cal F_{r, s}$ denote the closed subspace of $C^*(S_X)$ generated by
$\{S_xS_y^* ; x  \in X^{\times s}, y\in X^{\times r}\}$.
  Let, for $T \in \K = \Cal K_\A(X_\A) = \Cal L_\A(X_\A)$, $\quad \pi _\K(T) \in
\Cal L_\A(F(X)_\A)/\Cal K_\A(F(X)_\A)$ denote  the  image in the quotient algebra
of the operator on $F(X)$ that acts trivially on $X^{\otimes 0}=\A$ and that on 
$X^{\otimes r}$, $r = 1,2,\dots$  acts by
$
x_1\otimes\dots\otimes x_r \mapsto (Tx_1)\otimes\dots\otimes x_r.
$
  Since $\pi _\K(\theta _{x,y}) = S_xS_y^*$ for $x,y \in X$,
we have $\pi _\K(T) \in \Cal F_{1, 1} \subset C^*(S_X)$ for $T \in \K$.
  Put  $\pi_\A(a) = \pi_\K(\phi(a))$ for
$a \in \A$. Then $(S, \pi_\K, \pi_\A)$ is a representation of
$_\K X_\A$ in $C^*(S_X)$ {\it with support I\/}.
It has been shown in \cite{Pi} that $C^*(S_X)$ is the universal
$C^*$--algebra generated by such a representation,
 therefore it is isomorphic to   $\Cal O_X$.

  Let $\{u_1,\dots,u_n\}$ be a basis of $X$.
  Then the identity $$S_xS_y^* = \sum_{i=1}^n S_xS_{u_i}S_{u_i}^*S_y^*\ ,$$
for $x \in X^r$, $y\in X^s$, shows that $\Cal F_{r, s} \subset \Cal F_{r+1, s+1}$.
  We denote by ${\Cal F_{\infty}}^{(k)}$ the norm closure of $\cup _{r=0}^{\infty}
\Cal F_{r, r+k}$.
The unitary representation $U$ of $\Bbb T$ on $F(X)$ defined by
 $U(t)x=t^rx$, $t\in{\Bbb T}\ ,$ $x\in X^{\otimes r}$
implements an automorphic  action $\gamma : \Bbb T
\rightarrow Aut \ \Cal O _X$ such that $\gamma_t(S_x) = tS_x$
for $t \in \Bbb T$ and $x \in X$, thus making $\O_X$ into a 
${\Bbb Z}$--graded $C^*$--algebra.
  The fixed point algebra $\Cal O_X^{\gamma}$ is 
${\Cal F_{\infty}}^{(0)}$.
  We denote by $m_0$ the faithful conditional expectation of $\Cal O_X$
onto ${\Cal F_{\infty}}^{(0)}$ defined by
   $m_0(T) = \int _{\Bbb T} \gamma _t(T)dt$ for $T \in \Cal O_X$.

   We now recall  the  construction of  $\O_X$ by \cite{DPZ},
 intrinsic to the category of 
Hilbert $C^*$--bimodules. A general  construction  \cite{DR3}
 associates functorially to each object $\rho$ of a strict tensor 
$C^*$--category, a $C^*$--algebra $\O_\rho$.
This applies, with no substantial modifications, 
to the case where $\rho$ is just an abject 
of a strict {\it right  tensor} $C^*$--category $\T$. More specifically, $\T$
is a $C^*$--category
for which  the set of objects is a  semigroup with unit
$\iota$, such that  for any triple of objects $\rho\ , \sigma\ , \tau$
there are maps ${L_\tau}^{\rho, \sigma}$ defined on the space $(\rho, \sigma)$
of arrows from $\rho$ to $\sigma$ (``tensoring on the right by $1_\tau$'') 
that satisfy:
$$
{L_\tau}^{\rho, \sigma}: (\rho, \sigma)\to(\rho\tau, \sigma\tau)
$$
$${L_\omega}^{\rho\tau , \sigma\tau}{L_\tau}^{\rho, \sigma}={L_{\tau\omega}}^{\rho, \sigma}\ ,$$
$${L_\iota}^{\rho, \sigma}=id\ .$$
These maps are assumed to satisfy  natural compatibility
relations with the $^*$--involution and composition of the category $\T$.
Any strict tensor $C^*$--category is a strict right tensor $C^*$--category
by ${L_\tau}^{\rho, \sigma}(T)=T\otimes 1_\tau$. 
For simplicity, we assume that ${L_\tau}^{\rho, \sigma}$
are isometric mappings.

Let ${^0{\O_\rho}}^{(k)}$ be the algebraic inductive limit of linear spaces
$(\rho^r, \rho^{r+k})$ by the maps ${L_\rho}^{\rho^r, \rho^{r+k}}$.
Then the algebraic direct sum
 $^0{\O_\rho}=\oplus_{k\in{\Bbb Z}}{^0{\O_\rho}^{(k)}}$ is a
$\Bbb Z$--graded $^*$--algebra. Let
$\gamma: {\Bbb T}\to \Aut^0{\O_\rho}$ be the $\Bbb T$--action defining the grading.
The algebra  $\O_\rho$
is defined to be  the completion of $^0{\O_\rho}$ in the unique $C^*$--norm
for which $\gamma_z$ is
isometric.

Let $\A\subset\D$ be a unital inclusion of $C^*$--algebras.
A typical example of strict right   tensor $C^*$--category is the category
$\H_\A(\D)$ of   Hilbert $\A$--bimodules   contained in a
$C^*$--algebra $\D$, finite projective as right $\A$--modules and with support $I$.
 An object  of $\H_\A(\D)$ is
a closed subspace $X$ of
 $\D$ for which $\A X=X\A=X$ and $X^*X\subseteq \A$.  $X$ carries an obvious Hilbert
 $\A$--bimodule structure.  
The set of arrows  between two objects $X$ and $Y$
is the 
  space $\L_\A(X_\A, Y_\A)=\K_\A(X_\A, Y_\A)$ 
of right adjointable $\A$--module maps, identified with the 
subspace $YX^*$ of $\D$.
The operations of composition and $^*$--involution  on the arrows
of  $\H_\A(\D)$
 are given by the operator product and $^*$--involution inherited from $\D$.
For any pair  of objects $X$, $Y$ of $\H_\A(\D)$, $XY$, 
 the linear span of operator products
$xy$, $x\in X$, $y\in Y$, is a strict realization of
the tensor product $X\otimes_\A Y$ of Hilbert $\A$--bimodules. 
The unit object is the trivial bimodule
 $\A$. 
On arrows, the ``tensor product with the identity arrows''
is defined by
 ${L_Z}^{X, Y}(T)=T$.
Note that the Hilbertian norm on $X$ coincides 
with  the norm inherited from
$\D$.
 
 Doplicher--Roberts construction applied to an object $X$ of $\H_\A(\D)$
provides a $C^*$--algebra containing    $X$ and 
$\A$. Since $X_\A$ is finite projective,
 the subspace   $^0{{\O_\rho}^{(k)}}$ of $\O_\rho$ identifies 
with  $\cup_r X^{r+k}{X^r}^*\subset \D$,
therefore $\O_\rho$ is a $C^*$--algebra generated by $X$
(that in general
 fails to identify
with 
the $C^*$--subalgebra $C^*(X)$ of $\D$ generated by $X$). 
It has been shown in \cite{DPZ} that this $C^*$--algebra is isomorphic to $\O_X$.
For  convenience, we will give a proof of this fact.
  We  need the following standard computation of the norm of a finite rank operator.

\proclaim{2.1 Lemma}
   Let $\A$ be a $C^*$--algebra and $X_\A$ a right Hilbert $\A$--module.
   For $x_1,\dots,x_n$, $y_1,\dots,y_n \in X$, we have
$$
   \Vert \sum_{i=1}^n \theta_{x_i,y_i} \Vert
 = \Vert ((x_i\vert x_j)_\A)_{ij}^{1/2} ((y_i \vert y_j)_\A)_{ij}^{1/2} \Vert
$$
where the right hand side is the norm on $M_n(\A)$.
   In particular, if $X = \A$, then for $x_1,\dots,x_n, y_1,\dots,y_n \in \A$
we have
$$
   \Vert \sum_{i=1}^n x_iy_i^* \Vert
 = \Vert (x_i^*x_j)_{ij}^{1/2}(y_i^*y_j)_{ij}^{1/2} \Vert
$$
\endproclaim
\demo{Proof}
   For $x,y \in X$, we have
$$
\align
&\Vert \theta_{x,y} \Vert^2 = \Vert \theta_{x,y}^*\theta_{x,y} \Vert
= \Vert \theta_{y(x\vert x)_\A,y}\Vert=\\ 
&\Vert \theta_{y(x\vert x)_\A^{1/2},y(x\vert x)_\A^{1/2}} \Vert= 
\Vert (y(x\vert x)_\A^{1/2}\, \vert \, y(x\vert x)_\A^{1/2})_\A\Vert= \\
& \Vert(x\vert x)_\A^{1/2}(y\vert y)_\A^{1/2} \Vert^2.
\endalign
$$
   Consider $Y := X \oplus \cdots \oplus X = X \otimes \Bbb C^n$ as
a right $\A \otimes M_n(\Bbb C)$--module as in the proof of
[KW1; Proposition 1.18].
   For $x = (x_1,\ldots ,x_n), y = (y_1,\ldots ,y_n) \in Y$, we put
$(x\vert y)_{\A\otimes M_n(\Bbb C)} = ((x_i\vert y_j)_\A)_{ij}$.
   Let $$\Theta_{x,y}(z) = x(y\vert z)_{\A\otimes M_n(\Bbb C)}$$
for $x,y,z \in Y$.
   Then $\Theta_{x,y} = (\sum_{i=1}^n \theta_{x_i,y_i})\otimes I$
by the identification 
$$\Cal L_{\A\otimes M_n(\Bbb C)}(Y_{\A\otimes M_n(\Bbb C)})=
\Cal L_\A(X_\A) \otimes \Bbb C.$$
   Therefore we have
$$
   \Vert \sum_{i=1}^n \theta_{x_i,y_i} \Vert
 = \Vert \Theta_{x,y} \Vert
 = \Vert (x\vert x)_{\A\otimes M_n(\Bbb C)}^{1/2}
         (y\vert y)_{\A\otimes M_n(\Bbb C)}^{1/2} \Vert 
$$
$$
=\Vert ((x_i\vert x_j)_\A)_{ij}^{1/2}((y_i\vert y_j)_\A)_{ij}^{1/2} \Vert
$$
\enddemo

   The following Lemma  shows that any representation of a right Hilbert
$C^*$--module $X_\A$ in a $C^*$--algebra 
extends to a  representation of the  imprimitivity bimodule $_\K X_\A$
in the same $C^*$--algebra.

Let $\T$ denote a full $C^*$--subcategory of the category $\H_\A$ of  finite 
projective right Hilbert modules over $\A$. 

\proclaim{2.2. Lemma}
([Pi; Lemma 3.2]) \quad Let $\A$ and $\D$ be unital   
$C^*$--algebras.
 Let be given a linear map $\pi_X : X \rightarrow \D$
for each object $X$ of $\T$ and a  unital $^*$--homomorphism
$\pi_\A : \A \rightarrow \D$ 
  such that
$$
\pi_X(xa) = \pi_X(x)\pi_\A(a) \quad \text{and} \quad
\pi_X(x)^*\pi_X(y)=\pi_\A((x\vert y)_\A)  
$$
for $x,y \in X$ and $a \in \A$.
  Then there exists a unique $^*$--functor  $\pi : \T  \rightarrow \H_\A(\D)$ 
with full image that associates $\pi_X(X)$ with any object
 $X\in\T$, and such that
$$
  \pi(\theta_{y, x}) = \pi_Y(y)\pi_X(x)^* \quad 
\text{and} \quad \pi_X(kx) = \pi_\K(k)\pi_X(x)
$$
for $x\in X$, $y\in Y$ and $k \in \K_\A(X_\A, Y_\A)$.
Furthermore if $\pi_\A$ is one to one, then $\pi$ is faithful
and $\pi_X$ is an isometry.
\endproclaim
\demo{Proof}
For $x_1,\dots,x_n\in X, y_1,\dots,y_n\in Y$, Lemma 2.1 with a slight modification 
shows that
$$
\align
&\Vert \sum_{i=1}^n \theta_{x_i,y_i} \Vert
  = \Vert ((x_i\vert x_j)_\A)_{ij}^{1/2} ((y_i\vert y_j)_\A)_{ij}^{1/2} \Vert
  \\
&\geqq\Vert (\pi_\A((x_i\vert x_j)_\A))_{ij}^{1/2}
              (\pi_\A((y_i\vert y_j)_\A))_{ij}^{1/2} \Vert\\
&=\Vert (\pi_X(x_i)^*\pi_X(x_j))_{ij}^{1/2}
        (\pi_Y(y_i)^*\pi_Y(y_j))_{ij}^{1/2} \Vert\\
&=\Vert \sum_{i=1}^n \pi_X(x_i)\pi_Y(y_i)^* \Vert.
\endalign
$$
   The above estimate shows that there exists a contraction
$\pi : \K_\A(X_\A, Y_\A) \rightarrow \D$ such that
$\pi(\theta_{y,x}) = \pi_Y(y)\pi_X(x)^*$.
It is easily seen that the map $\pi: \T\to\H_\A(\D)$   is in fact a $^*$--functor.
Since
$$
   \pi_X(\theta_{x,y}(z)) = \pi_X(x(y\vert z)_\A)
=  \pi_X(x)\pi_\A((y\vert z)_\A)  
$$
$$
=\pi_X(x)\pi_Y(y)^*\pi_X(z)
=  \pi(\theta_{x,y})\pi_X(z)
$$
   we have that $\pi_X(kz) = \pi(k)\pi_X(z)$ for $k \in
 \K_\A(X_\A, Y_\A)$ and
$z \in X$.
   Furthermore assume that $\pi_\A$ is one to one.
   Then the only inequality in the above calculation becomes  an
equality. 
   Thus $\pi : \K_\A(X_\A, Y_\A) \rightarrow \D$ is also one to one.
   Since
$$
\Vert \pi_X(x) \Vert ^2 = \Vert \pi_X(x)^*\pi_X(x) \Vert
= \Vert \pi_\A((x\vert x)_\A)\Vert = \Vert (x\vert x)_\A \Vert
= \Vert x \Vert^2,
$$
$\pi$ is an isometry.
\enddemo

Let $\T_X$ denote the full subcategory of the category $\H_\A$ of finite
projective right Hilbert $\A$--modules, with objects 
$\{X^{\otimes r}, r=0,1,2,\dots\}$.

\proclaim{2.3. Corollary}
(\cite{Pi}) Let $\A$ be a unital
$C^*$--algebra and $_\A X_\A$ a 
Hilbert  $\A$--bimodule, full and finite projective as a right $\A$--module.
   Let $\phi : \A \rightarrow \K = \Cal K_\A(X_\A) = \Cal L_\A(X_\A)$ be the defining left
action of $\A$.
   Let $\Cal O_X = C^* \{ S_x ; x \in X \} $ and let $\pi_\A : \A \rightarrow
\Cal O_X$ be the canonical embedding.
   Then there exist representations $(\pi_r, \pi_\A)$, $r=0, 1, 2, \dots,$
 of Hilbert $\A$--bimodules
$X^{\otimes r}$ in $\O_X$ with support $I$ such that
$$
  \pi_r(x_1\otimes \cdots \otimes x_r) = S_{x_1}\cdots S_{x_r}.
$$
Furthermore there is a faithful right tensor 
$^*$--functor $\pi: \T_X\to\H_\A(\O_X)$
such that
$\pi(X^{\otimes r})=\pi_r(X^{\otimes r})$ and
$$ 
\pi:\theta_{x_1\otimes \cdots \otimes x_s, y_1,\otimes \cdots \otimes y_r}\in
\Cal K_\A({X^{\otimes r}}_\A,
 {X^{\otimes s}}_\A) \rightarrow  S_{x_1}\cdots S_{x_s}S_{y_r}^*\cdots S_{y_1}^*\in
\Cal F_{r, s}.
$$
\endproclaim
\demo{Proof}
   Apply Lemma 2.2 to ${X^{\otimes r}}$ in place of
$X$ and put $\D = \Cal O_X$.
   Since
$$
  S_{x_r}^*\cdots S_{x_1}^*S_{y_1}\cdots S_{y_r}
= \pi_\A((x_1\otimes \cdots \otimes x_r \ \vert \
         y_1\otimes \cdots \otimes y_r )_\A),
$$
there exists an isometry $\pi_r : {X^{\otimes r}} \rightarrow  \Cal O_X$ such that
$\pi_r(x_1\otimes \cdots \otimes x_r)
= S_{x_1} \ldots S_{x_r}$ that
 satisfies 
$$
 \pi_r(xa) = \pi_r(x)\pi_\A(a),\quad
 \pi_r(\phi(a)x)=\pi_\A(a)\pi_r(x)
$$
and
$$
 \pi_r(x)^*\pi_r(y) = \pi_\A((x \vert y)_\A),
$$
hence $\pi_r$ is a faithful representation of $X^{\otimes r}$ in $\O_X$, clearly with support
$I$.
By Lemma 2.2, there exists  a faithful $^*$--functor $\pi: \T_X\to\H_\A(\D)$
 that satisfies the stated properties. We show that $\pi$ is right tensor.
Let $\{u_1, \dots, u_n\}$ be a basis of $X$.
$$
\align
 & \quad \pi (\theta_{x_1\otimes \cdots \otimes x_s, y_1\otimes \cdots \otimes 
y_r} \otimes I)
 = \pi (\sum_{i=1}^n \theta_{x_1\otimes \cdots \otimes x_s \otimes u_i,
                             y_1\otimes \cdots\otimes y_r \otimes u_i}) \\
 & =\sum_{i=1}^n S_{x_1}\ldots S_{x_s}S_{u_i}
                  S_{u_i}^*S_{y_r}^*\ldots S_{y_1}^*=S_{x_1}\ldots S_{x_s}
S_{y_r}^*\ldots S_{y_1}^*.
\endalign
$$
\enddemo

Let $\rho=X_\A$ be a Hilbert 
$\A$--bimodule contained in a $C^*$--algebra
$\D$ with support $I$,
 regarded as an object of the strict right tensor $C^*$--category
$\H_\A(\D)$, and associate to it the Doplicher--Roberts algebra $\O_{\rho}$. 

\proclaim{2.4. Lemma} 
Let ${\O_\rho}^{(k)}$ denote the closure of $^0{\O_\rho}^{(k)}$ in $\O_{\rho}$.
There is an isomorphism of ${\Bbb Z}$--graded $^*$--algebras
$$
\psi: 
\oplus_{k\in{\Bbb Z}}{\O_\rho}^{(k)}\to\oplus_{k\in{\Bbb Z}}{\Cal F_{\infty}}^{(k)}
$$
such that  $\psi(x)=S_x, x\in X.$ 
Furthermore $\psi$ is compatible with the inclusions
$\Cal K_\A({X^{\otimes r}}_\A, {X^{\otimes r+k}}_\A)
\hookrightarrow \Cal K_\A({X^{\otimes r+1}}_\A, {X^{\otimes r+k+1}}_\A)$ in 
$\O_{\rho}^{(k)}$ and
$\Cal F_{r,r+k} \subset \Cal F_{r+1,r+k+1}$ in ${{\Cal F}_{\infty}}^{(k)}$.
\endproclaim
\demo{Proof}
   Since the $^*$--functor $\pi$ defined in  Corollary 2.3 is
 compatible with right tensorization with the identity arrows, we may define
a $^*$--isomorphism of ${\Bbb Z}$--graded $^*$--algebras
$
\psi: 
\oplus_{k\in{\Bbb Z}}{\O_\rho}^{(k)}\to\oplus_{k\in{\Bbb Z}}{\Cal F_{\infty}}^{(k)}
$
by letting $$\psi|_{\K({X^{\otimes r}}_\A, {X^{\otimes s}}_\A)}=\pi\ .$$
\enddemo
 
As a consequence of the previous lemma we deduce the following result,
proved also in \cite{DPZ}.

\proclaim{2.5. Proposition}
The isomorphism 
$\psi: \oplus_{k\in{\Bbb Z}}{\O_\rho}^{(k)}\to\oplus_{k\in{\Bbb 
Z}}{\Cal F_{\infty}}^{(k)}$ 
of ${\Bbb Z}$--graded $^*$--algebras defined in the previous lemma extends
to an isomorphism $\psi: \O_\rho\to\O_X$ of $C^*$--algebras.
\endproclaim
\demo{Proof} In view of Lemma 2.4 we need
 to show that the norm of $\oplus_{k\in{\Bbb 
Z}}{\Cal F_{\infty}}^{(k)}$  inherited from 
$\Cal L_\A(F(X)_\A)/\Cal K_\A(F(X)_\A)$
makes the canonical action of ${\Bbb T}$ continuous. But this is clear since this 
action is  implemented by the quotient image
of a unitary representation of $\Bbb T$
on  $F(X)_\A$.
\enddemo

In the following we shall identify $\O_X$ with the algebra $\O_\rho$,
omitting the isomorphism $\psi$.
In particular, $X$ and $\A$, as well as the spaces
$\K_\A({X^{\otimes r}}_\A, {X^{\otimes r}}_\A)$,  will be regarded as closed subspaces 
of $\O_X$.

   The existence of the following sequence of conditional expectations
will be useful later.

\proclaim{2.6. Lemma}
  Consider the same situation as in Corollary 2.3. 
  Assume furthermore  that $X$ is a Hilbert $\A$--bimodule of finite type
in the sense of \cite{KW1} with  left inner product ${}_\A(\ \vert \ )$.
  Then there exist conditional expectations
$E_r: {\O_X}^{(0)} \rightarrow \Cal K_\A({X^{\otimes r}}_\A)$ \  (for $r = 0,1,2,\ldots$)
such that $T = \lim_{r\to\infty}E_r(T)$ (in the norm topology) for $T \in
{\O_X}^{(0)}$.
   Moreover there exists a conditional expectation
$E^{\Cal O_X}_\A : \Cal O_X \rightarrow \A$
\endproclaim  
\demo{Proof}
As in  [KW1, Lemma 3.23, 3.24],  there exists
a conditional expectation
$E_r^{r+k} : \Cal K_\A({X^{\otimes r}}_\A \otimes _\A {X^{\otimes k}}_\A) \rightarrow
\Cal K_\A({X^{\otimes r}}_\A)$
such that 
$$E_r^{r+k}(\theta_{x_1\otimes y_1, x_2\otimes y_2})= 
(r-Ind[X^{\otimes k}])^{-1} \theta_{x_1{}_\A(y_1\vert y_2),x_2}$$
 \ for \ $x_1, x_2 \in {X^{\otimes
r}}$
\ and \ $y_1, y_2 \in {X^{\otimes k}}$.
   We also have a conditional expectation $G : \Cal K_\A(X_\A) \rightarrow
\A$ such that $G(\theta_{x,y}) = (r-Ind[X])^{-1} {}_\A(x\vert y)$ for
$x,y \in X$ by  [KW1; Lemma 1.26]. 
   Then we can easily construct the desired conditional expectations
$E_r : {\O_X}^{(0)} \rightarrow \Cal K_\A({X^{\otimes r}}_\A)$.
  Put $E_0 = G\circ E_1$ and $E^{\Cal O_X}_\A = E_0\circ m_0$, where
$m_0$ is the conditional expectation of $\Cal O_X$ onto ${\O_X}^{(0)}$ obtained
overaging over the action of ${\Bbb T}$.
\enddemo 

     Lemma 2.6 shows that, if $X$ is a  Hilbert $\A$--bimodule
 of finite type,  the relative commutant algebra
$\A'\cap {\O_X}^{(0)}$ is approximated by the subalgebras
$\cup _{r=1}^{\infty} \A'\cap \K_\A({X^{\otimes r}}_\A)$ 
(cf \cite{KW1}).
 
\heading
3. Canonical Completely Positive Maps
\endheading

   Let ${u_1,\ldots ,u_n}$ be a finite basis of $X = X_\A$.
     We define a completely positive map $\sigma : \Cal O_X \rightarrow \Cal O_X$
by $\sigma (T) = \sum_{i=1}^n u_iT{u_i}^*$.
   We note that $\sigma$ does depend on the choice of the basis,
and that it is not a $^*$--endomorphim in general. However,

\proclaim{3.1. Lemma}
   The restriction of $\sigma$ to $\A'\cap \Cal O_X$ is a unital
 $^*$--monomorphism that does not depend on the choice of the basis
of $X_\A$.
\endproclaim
\demo{Proof}
  For $T_1, T_2 \in \A'\cap {\Cal O}_X$,
$$
\sigma (T_1)\sigma (T_2)
   = \sum_i u_iT_1{u_i}^* \sum_j u_jT_2{u_j}^* 
$$
$$ 
=\sum_i \sum_j {u_i}T_1(u_i\vert u_j)_\A T_2{u_j}^*
   = \sum_i \sum_j {u_i}(u_i\vert u_j)_\A T_1T_2{u_j}^* 
$$
$$
 =\sum_j {u_j}T_1T_2{u_j}^* = \sigma (T_1T_2).
$$
   Suppose that $\sigma (T) = 0$ for $T \in \A'\cap \Cal O_X$.
   Then $\sigma (TT^*) = \sum_i {u_i}TT^*{u_i}^* = 0$.
   Hence ${u_i}T = 0$.
   Then ${u_ia}T = {u_i}aT = {u_i}Ta = 0$ for $a \in \A$.
   Thus $xT = 0$ for all $x \in X$.
   Therefore $(y\vert x)_\A T = y^*xT = 0$.
   Since $X_\A$ is full, $T = 0$.
   Thus $\sigma$ is isometric on $\A'\cap \Cal O_X$.
   Let $\{ v_1,\ldots ,v_m \}$ be another basis of $X_\A$.
   Then ${v_j} = {\sum_i u_i(u_i\vert v_j)_\A}$.
   For $T \in \A'\cap \Cal O_X$, we have
$$
\sum_{j=1}^m {v_j}T{v_j}^*
 = \sum_j (\sum_i {u_i}(u_i\vert v_j)_\A) T
          (\sum_k {u_k}(u_k\vert v_j)_\A^*) 
$$
$$
=\sum_k\sum_i {u_i}(\sum_j(u_i\vert v_j)_\A(v_j\vert u_k)_\A
     T{u_k}^*
   = \sum_k \sum_i {u_i}(u_i\vert u_k)_\A T{u_k}^* 
$$
$$
=\sum_k {u_k}T{u_k}^* = \sigma (T).
$$
\enddemo

\proclaim{3.2. Lemma}
   For $T \in\A'\cap \Cal O_X$ and $x_1,\ldots x_m \in X$ we have
$$\sigma ^m(T){x_1}\ldots {x_m} = {x_1}\ldots {x_m}T$$ and
$\sigma ^m(T)$ commutes with $\K_\A({X^{\otimes m}}_\A)$.
   In particular $\sigma (T)$ commutes with $\A  \subset \K_\A({X}_\A)$
and $\sigma$ preserves $\A' \cap \Cal O_X$.
\endproclaim
\demo{Proof}
$$
\align
 & \quad \quad  \sigma ^m(T){x_1}\ldots {x_m}
= \sum_{i_m} \ldots \sum_{i_1}
u_{i_m}\ldots u_{i_1}Tu_{i_1}^*\ldots u_{i_m}^*{x_1}\ldots {x_m} \\
& = \sum_{i_m}\ldots \sum_{i_1}u_{i_m}\ldots u_{i_1}u_{i_1}^*\ldots
u_{i_m}^*{x_1}\ldots {x_m}T ={x_1}\ldots {x_m}T.
\endalign
$$
Therefore for $x_1,\ldots,x_m,y_1,\ldots,y_m\in X$,
$$
  \sigma ^m(T){x_1}\ldots {x_m}{y_m}^*\ldots {y_1}^*
= {x_1}\ldots {x_m}T{y_m}^*\ldots {y_1}^*
$$
$$
={x_1}\ldots {x_m}{y_m}^*\ldots {y_1}^*\sigma ^m(T)
$$
\enddemo

   The following Lemma shows that $\sigma$ acts as a shift operator on the subpaces
$\L_\A({X^{\otimes r}}_\A, {X^{\otimes s}}_\A)$.

\proclaim{3.3. Lemma}
Let 
$$
\psi: 
\oplus_{k\in{\Bbb Z}}{\O_\rho}^{(k)}\to\oplus_{k\in{\Bbb Z}}{\Cal F_{\infty}}^{(k)}
$$
 be the
$^*$--isomorphism defined
in Lemma 2.4.  For $T \in \A'\cap \Cal K_\A({X^{\otimes m}}_\A)$, we have
$\psi ^{-1}\sigma \psi (T) = I \otimes T \in \Cal K_\A(X\otimes _\A{X^{\otimes m}}_\A)$
\endproclaim
\demo{Proof}
   Since $T$ commutes with the left action of $\A$, $I\otimes T$ makes sense.
   Let $T = \sum_{x,y} \theta _{x,y}$, where $x = x_1 \otimes \cdots \otimes
x_m$ \ and \  $y = y_1 \otimes \cdots \otimes y_m \in {X^{\otimes m}}$ run over  a finite
set of simple tensors.
   Then
$$
\align
  &  \quad \quad \psi \sigma \psi ^{-1}(T)
= \psi \sigma \psi ^{-1}(\sum_{x,y} \theta _{x,y})
= \psi \sigma (\sum_{x,y} S_xS_y^*) \\
  & = \psi (\sum_{i=1}^n \sum_{x,y}S_{u_i}S_xS_y^*S_{u_i}^*)
= \sum_{i=1}^n \sum_{x,y} \theta_{u_i\otimes x,u_i\otimes y}
= I \otimes T.
\endalign$$
\enddemo
  
\heading
4. The Ideal Structure of $\Cal O_X$
\endheading

   In this section we study the ideal structure of $\O_X$ under certain 
assumptions on the bimodule $X$. In particular, we
get a simplicity criterion for $\Cal O_X$.

   It was shown in \cite{El}, \cite{Ki} that
  discrete crossed product $C^*$--algebras are simple if
the group action satisfies an outerness property. On the other hand,
Cuntz--Krieger algebras
$\O_A$ are simple if the defining matrix $A$ 
is irreducible and non permutation \cite{CK}. More generally, if $A$ satisfies
property $(I)$ of \cite{CK} then $\O_A$ is uniquely determined by generators 
and relations.
Furthermore if $\A$ satisfies property $(II)$  (which is stronger than
$(I)$) then the ideal structure of $\O_A$ is explicitly determined in  \cite{Cu2}.  
    Matsumoto   generalized the simplicity argument
 to  $C^*$--algebras
associated with subshifts \cite{Ma}.
Furthermore,  Doplicher--Roberts algebras $\O_\rho$ were shown to be
 simple in the case where   $\rho$ is
a unitary finite dimensional representation of a compact group
 of determinant 1. 

 Our aim is to formulate a criterion to study the ideal structure of $\O_X$ that
generalizes the mentioned results.

\proclaim{4.1. Definition}  Let $\A$ be a unital $C^*$--algebra and $X_\A$ a finite projective
right Hilbert $\A$--module.
   Let $\phi : \A \rightarrow \K = \Cal K_\A(X_\A) = \Cal L_\A(X_\A)$ be a unital
isometric $^*$--homomorphism.
   Then a closed ideal $J$ of $\A$ is called $X$--invariant if
$(x\vert \phi(a)y)_\A \in J$ for $x, y\in X$ and $a\in J$.
\endproclaim

 We note that
$J_X:=\{a\in \A: (x\vert \phi(a)y)_\A\in J, x, y\in X\}$ is a closed $X$--invariant
 ideal of $\A$
containing $J$, if $J$ is $X$--invariant. Note that if $X$ is not full and
$J$ is the proper ideal generated by inner products then
$J_X=\A$. Otherwise, if $X$ is full then $J_X$ is proper if $J$ is. In particular,
if $J$ is a maximal proper $X$--invariant ideal then $J_X=J$.
   An algebra $\A$ with  no proper $X$--invariant ideal will be called
{\it $X$--simple}.

For any closed  ideal $J$ of $\A$  we denote by $X_J$ 
the closed subset of $X$
of elements $x\in X$ for which $(x\vert x)_\A\in J$. If we approximate
$x$ with $xu_\alpha$, with $u_\alpha$ an approximate unit of $J$, we see that
$(y\vert x)_\A\in J$ for all $y\in X$. This shows in particular that $x+y\in X_J$ 
if $x, y\in X_J$, i.e. $X_J$ is a subspace of $X$.
Furthermore,  since
 $XJ\subseteq X_J$ and $X_J\A\subseteq X_J$,   $X/X_J$ is a right Hilbert 
module over $\A/J$ in a natural way.
 
\proclaim{4.2. Proposition} If $J$ is a closed $X$--invariant ideal of $\A$ then 
$\tilde{\phi}:  \A/J\to\L_{\A/J}(X/X_J)$, $\tilde{\phi}([a])[x]=[\phi(a)x]$
defines a $^*$--homomorphism such that
 $ker \tilde{\phi}=J_X/J$.
In particular, if  $J_X=J$ 
then $\tilde{\phi}$
is isometric.
\endproclaim
\demo{Proof}
Since $X$ is finite projective and $J$ is $X$--invariant, 
$\phi(J)X\subseteq XJ\subseteq X_J$.
Furthermore for $x\in X_J$, $a\in \A$, 
$(\phi(a)x\vert \phi(a)x)=(x\vert \phi(a^*a)x)\in J$, thus $\phi(\A)X_J\subseteq X_J\ .$
 It follows
that $\tilde{\phi}$ is a well defined $^*$--homomorphism. The rest of the proof
is now clear.
\enddemo

   A bimodule $X$ is called $(I)$--{\it free} if    for any $k \in \Bbb N$
there exists  an element
$T_k \in \A'\cap \L_\A({X^{\otimes p_k}}_\A, {X^{\otimes q_k}}_\A)$ 
of norm $1$ satisfying the following
conditions: 
 $$\A \owns a \mapsto \phi(a){T_k}^*T_k \in \K_\A({X^{\otimes r_k}}_\A)\eqno(4.1)$$ 
is completely isometric and 
$$\|{T_k}^*\sigma^k(T_k)\|<1\ .\eqno(4.2)$$

A bimodule $X$ is called $(II)$--{\it free} if for any closed $X$--invariant ideal $J$ of $\A$
such that $J_X=J$, $X/X_J$ is $(I)$--free. If $J=\{0\}$ then $J_X=\{0\}$ since
$\phi$ is faithful, therefore $(II)$--freeness implies $(I)$--freeness.

Let $(V, \rho_\A)$ be a representation of the Hilbert $\A$--bimodule $X$
in $\D$ with suppport $I$.
Let $C^*(V(X))$ and $^*alg(V(X))$ denote respectively
the $C^*$--subalgebra and the $^*$--subalgebra of $\D$ generated by $V(X)$.

   The following Theorem   shows that if $X$ is $(I)$--free, then 
generators and relations determine
$\Cal O_X$ uniquely,  whenever $\A$ is  faithfully represented. If in addition $X$ is 
$(II)$--free then one can determine explicitly the ideal structure of $\O_X$.

\proclaim{4.3. Theorem} Let $X$ be a full Hilbert module over a $C^*$--algebra
$\A$ and
   let $\D$ be a unital $C^*$--algebra, $(V, \rho_\A)$ a representation
of $_\A X_\A$ in $\D$ with support $I$ and
   let $\varphi : \Cal O_X \rightarrow \D$ be the unique $^*$--homomorphism
such that $\varphi(x) = V_x$ and  $\varphi(a) = \rho_\A(a)$
for $x \in X, a \in \A$
(by the universality of $\Cal O_X$).
\par \noindent
 $(i)$ If $\A \owns a \rightarrow \rho_\A(a) \in \D$ is one to one, then
the restriction of $\varphi$ to ${\O_X}^{(0)}$ is one to one.
Furthermore, if  $X$ is $(I)$--free
   then 
\roster
\item $\varphi : \O_X\to \D$ is one to one,
hence $^0\O_X\cong$ $^*$--alg$(V(X))$,
\item $^0\O_X$ has a unique $C^*$--norm.
\endroster
$(ii)$ If $X$ is $(II)$--free then
\roster 
\item' the maps $\J \mapsto J:=\J\cap\A$, 
$J\to\J:=cls\{X^rJ{X^s}^*, r,s=0,1,2,\dots\}$ are
 inclusion preserving
 bijective correspondences, inverse of one another, between closed ideals of $\O_X$
and $X$--invariant ideals $J$ of $\A$ for which $J_X=J$.
\item' if $J$ is the ideal of $\A$ corresponding to the ideal $\J$ of
$\O_X$ as in $(1)'$ then  $\O_X/\J=\O_{X/X_J}$.
 \endroster
\endproclaim
\demo{Proof}
$(i)$:  Apply Lemma 2.2 to ${X^{\otimes r}}$ in place of $X$, for any $r$.
 The restriction of $\varphi$ to $\K_\A({X^{\otimes m}}_\A)$ is isometric,
  hence the restriction of $\varphi$ to ${\O_X}^{(0)}$ is isometric.
  Let us assume  that $X$ is $(I)$--free.

We claim, without proving it, that the arguments of Lemma 2.2 can be easily
generalized to show that
for any pair of Hilbert bimodules $X$, $Y$ and any
$x_1,\dots,x_n$, $y_1,\dots,y_n \in X$, $T\in \phi(\A)'\cap\L_\A(Y_\A)$,
$$
   \Vert (\sum_{i=1}^n \theta_{x_i,y_i}\otimes 1_Y)\circ (1_X\otimes T) \Vert
 = \Vert ((x_i\vert x_j)_\A)^{1/2} ((y_i \vert y_j)_\A)^{1/2}diag((T^*T)^{1/2}) \Vert
$$
where the right hand side is the norm on $M_n(\phi(\A)'\cap\L_\A(Y_\A))$ and $diag(B)$
 is the matrix of $M_n(\A)$ with  $B$ in the diagonal.
It follows that  $(4.1)$ can be also read
$$\|\sigma^p({T_n}^*T_n)T\|=\|T\|\ ,\quad T\in \Cal K_\A({X^{\otimes p}}_\A)$$.

Let $r_n$ be sufficiently large so that ${T_n}^*\sigma^j(T_n)\in
\L_\A({X^{\otimes r_n}}_\A)$
for $j=1,\dots,n$.
We claim also that replacing $r_n$ by some other ${r'}_n$,
if necessary,  we may find ${T'}_n$ satisfying  the above property 
and the following property, stronger than $(4.2)$,
$$\|{{T'}_n}^*{\sigma}^n({T'}_n)\|<\varepsilon\ .$$
To prove the claim it is enough to assume $n=1$, since we
may  replace $T=T_1$ by $T_n$, $r_1$ by $r_n$ and
$X$ by $X^{\otimes n}$.
  Choose a natural number $q$ such that 
   $\|T^* \sigma (T) {\|}^{q+1} < \varepsilon \ .$
Then $ T' := T{\sigma}^{r_1}(T)\dots{\sigma}^{qr_1}(T)$ satisfies
the desired properties. Indeed, using the fact that
$ T^* \sigma (T) \in \L_\A({X^{\otimes r_1}}_\A)$ commutes with 
${\sigma}^{r_1}(T^*)$, ${\sigma}^{2r_1}(T^*)$, ... ,
${\sigma}^{qr_1}(T^*)$, we have that
$$    \| {T'}^* \sigma (T') \|
    = \| T^* \sigma (T){\sigma}^{r_1}(T^* \sigma (T))
      {\sigma}^{2r_1}(T^* \sigma (T)) \dots 
      {\sigma}^{qr_1}(T^* \sigma (T)) \|
$$
$$
\leq\| T^*\sigma(T) {\|}^{q+1} < \varepsilon \ .
$$
Moreover, for any $Z \in \Cal L_\A(X^{\otimes p}_\A)$, 
$$
   \| {\sigma}^p ({T'}^*T')Z \|
 = \| Z{\sigma}^p(T^*T){\sigma}^{p+r_1}(T^*T) \dots 
   {\sigma}^{p+qr_1}(T^*T) \| = \| Z \| \ .
$$

Let now $B=\sum_{j=-n}^n B_j$ be an element of $^0\O_X$, with
$$B_j\in\L_\A({X^{\otimes p}}_\A, {X^{\otimes p+j}}_\A)$$ and $n\neq0$. Note that
setting $$B':={\sigma}^p({{T'}_n}^*)B{\sigma}^p({{T'}_n}),$$
$$B'=
\sum_{j=-{n}}^{n}{B'}_j\ ,$$
where for non negative $j$,
$${B'}_j={\sigma}^p({{T'}_n}^*{\sigma}^j({{T'}_n})
{B}_j\ .$$
It follows that for some $p'$ and all $j,$
$${B'}_j\in\L_\A({X^{\otimes p'}}_\A, {X^{\otimes p'+j}}_\A).$$
Using the fact that the homogeneous subspaces ${^0\O_X}^{(k)}$ have a unique $C^*$--norm
and that $\varphi$ is faithful, hence isometric, on them, we deduce, by 
  property  $(4.1)$  that
$$\|\varphi({B'}_0)\|=\|{B'}_0\|=\|B_0\|,$$
furthermore,
$$\|{B'}_j\|<
\varepsilon\|B_j\|\quad\hbox{\rm for } j=-n, n,$$
$$\|{B'}_j\|\leq\|B_j\|\quad\hbox{\rm for } j=-(n-1),\dots, n-1,$$
$$\|\varphi(B')\|\leq\|\varphi(B)\|.$$
Applying   the same computation to $B'$, replacing ${T'}_n$ by
${T'}_{n-1},$ we deduce, after a finite number of steps, 
by the arbitrarity of $\varepsilon$, that
$$ \|B_0\|\ = \|\varphi(B_0)\| \leq\|\varphi(B)\|.$$ 
Thus there exists a conditional expectation $\widehat E : \varphi(\Cal O_X)
\rightarrow \varphi({\Cal O_X}^{(0)})$ such that 
$\widehat E (\varphi (B)) = \varphi (B_0)$ .  
Since $\widehat E \varphi = \varphi m_0$ ,  $\varphi$ is one to one on 
${\Cal O_X}^{(0)}$ and $m_0$ is faithful, 
we have that $\varphi$ is one to one,
and $(1)$ follows. On the other hand we have just shown that
the inequality $\|B_0\|\leq\|B\|$ holds
  in any $C^*$--norm on $^0\O_X$, hence the proof of $(2)$ is complete for,
 by \cite{DR1}, $^0\O_X$ has a unique $C^*$--norm with this property.

$(ii)$: 
Let now assume that $X$ is $(II)$--free and let $\J$ be a closed ideal
of $\O_X$. Then $J=\J\cap\A$ is a closed $X$--invariant ideal
of $\A$ satisfying $J_X=J$. Hence $X/X_J$ is $(I)$--free. Since $X_J = \J \cap X$,
 the bimodule $X/X_J$
is isometrically isomorphic to $X/\J$ thought of as a  Hilbert bimodule over
$\A/J$, thus there is a   natural representation
of $X/X_J$ in $\O_X/\J$  that is faithful on $\A/J$. We may apply a) to deduce that
$\O_{X/X_J}$ is isomorphic to $\O_X/\J$. 
Let
 $m_k$ denote the projection onto ${^0\O_X}^{(k)}$.
The first part of the proof applied to this
representation  also shows that if $B$ is an
element of  $\J$  then $m_0(B)$ is still in $\J$.
 Now
 for positive $k$, for any $x\in X^k$,
 $x^*m_k(B)=m_0(x^*B)$, thus  $x^*m_k(B)\in \J$,
hence $m_k(B)\in \J$. Since, by Fourier analysis,
 any element $B$ of $\O_X$ is the limit in norm of Cesaro sums of $m_k(B)$'s,
we deduce that any closed ideal $\J$
is invariant under the canonical action of ${\Bbb T}$, and hence it is 
generated by the subspaces $\J\cap\L_\A({X^{\otimes r}}_\A, {X^{\otimes s}}_\A)=
X^sJ{X^r}^*$, therefore the correspondence
is one to one. It remains to show that the correspondence is surjective, i.e.
that for any $X$--invariant ideal $J$ of $\A$ satisfying $J_X=J$ we have
$\J\cap\A=J$, where $\J$ is the closed ideal of $\O_X$ generated
by $J$. Clearly
$\J$ is invariant under the automorphic action of ${\Bbb T}$,
hence $J\subseteq \J^{(0)}\cap\A$. Now
$\J^{(0)}$ is the inductive limit of the $X^rJ{X^r}^*$, and 
$X^rJ{X^r}^*$ is a closed ideal of $\L_\A({X^{\otimes r}}_\A)$ 
that we will denote by $\J_r$. 
Let $\pi_r$ be the composition
of the inclusion of $\A$ in $\L_\A({X^{\otimes r}}_\A)$ with the quotient
map of $\L_\A({X^{\otimes r}}_\A)$ in $\L_\A({X^{\otimes r}}_\A)/\J_r$. 
Similarly, let $\pi_\infty$
be the composition of the inclusion of $\A$ in ${\O_X}^{(0)}$  with
the quotient map of ${\O_X}^{(0)}$ in ${\O_X}^{(0)}/\J^{(0)}$. 
 Then
ker$\pi_r=\A\cap\J_r$ and ker$\pi_\infty=\A\cap\J^{(0)}$, hence
there are obvious $^*$--monomorphisms
$$\A/\A\cap\J_r\to\L_\A({X^{\otimes r}}_\A)/\J_r$$
$$\A/\A\cap\J^{(0)}\to{\O_X}^{(0)}/\J^{(0)}.$$ 

If $b$ is an element of a $C^*$--algebra 
$\B$ and $\I$ is a closed ideal of $\B$, we denote by $[b]_{\B/\I}$ the 
image  of $b$ under the quotient map $\B\to \B/\I$. 

 Since a monomorphism 
between $C^*$--algebras is isometric, for any $a\in \A$
$$\|[a]_{\A/\A\cap\J^{(0)}}\|=\|[a]_{{\O_X}^{(0)}/\J^{(0)}}\|=
\hbox{\rm dist}(a, \J^{(0)})=$$
$$\lim_r \hbox{\rm dist}(a, \J_r)=
\lim_r\|[a]_{\L_\A({X^{\otimes r}}_\A)/\J_r}\|=
\lim_r\|[a]_{\A/\A\cap\J_r)}\|.$$
We claim that $J_X=J$ implies that $\A\cap\J_r=J$ for all $r$, 
hence adding this fact to the above chain
of equalities we deduce that 
$$\|[a]_{\A/\A\cap\J^{(0)}}\|=\|[a]_{\A/J}\|.$$ It follows
that the natural $^*$--homomorphism $\A/J\to \A/\A\cap\J^{(0)}$
is injective, i.e. $\A\cap\J^{(0)}=J$.
We finally prove the claim. If $a\in X^rJ{X^r}^*\cap\A$ 
then for $x, y\in X^r$, 
$x^*ay\in J$ hence for $x', y'\in X^{r-1}$,
 ${x'}^*ay'\in J_X=J$.
Iterating the argument we deduce that $a\in J$, and the proof is complete.
\enddemo

\proclaim{4.4. Corollary}
   Let $\A$ be a unital $C^*$--algebra and $_\A X_\A$ a  Hilbert $\A$--bimodule,
full and finite projective as a right Hilbert $\A$--module.
   Let $\phi : \A \rightarrow \K = \Cal K_\A(X_\A) = \Cal L_\A(X_\A)$ be a 
the defining left action of $\A$ on $X$.
   If $X$ is $(I)$--free and $\A$ is $X$--simple, then the $C^*$--algebra $\Cal O_X$
 is simple.
\endproclaim
\demo{Proof} Let $\J\not=\O _X$ be a closed ideal of $\O_X$.  Since
           $J:=\J \cap \A$ is $X$--invariant and $\A$ is $X$--simple, we have that
           $J=0$.  Then the  image of $\A$ in $\O_X/\J$ under the quotient map is faithful
           and the all the  relations are preserved.  Thus the
           quotient map is one to one by Theorem 4.3 $(1)$, that is  $\J = 0$.
           Therefore  $\O _X$ is simple. 
\enddemo

{\bf Remark 1.} Let $X$ be any Hilbert bimodule over $\A$.
 It has been proved in \cite{DPZ} that $X$--simplicity of $\A$
 is a necessary condition  to get simplicity of $\O_X$, but it is not
sufficient in general. Similarly to the case of
discrete crossed products \cite{El,Ki},
one need to assure furthermore that a certain
Connes spectrum naturally associated with $X$ is full \cite{DPZ}. 

{\bf Remark 2.}
   Following more closely the techniques of
 Elliott \cite{El} and Kishimoto \cite{Ki}, it is easy to see that
we can replace the assumption of $(I)$--freeness in Corollary 4.4 by the following
condition (i) or a  weaker condition (ii).

   (i) For any $k \in \Bbb N$, any $z \in \L_\A({X^{\otimes k}}_\A)$, any finite
subset $\{ z_1, \ldots , z_m \} \subset \L_\A({X^{\otimes k}}_\A)$ and  any
$\epsilon > 0$, there exists $q \in {\O_X}^{(0)}$ with
$\Vert q \Vert = 1$ such that
$\Vert qzq \Vert \geqq \Vert z \Vert - \epsilon $, \ 
$ \Vert qz_i - z_iq \Vert \leqq \epsilon $ for $i = 1, \ldots , m$
\quad and \quad $\Vert q\sigma ^j(q) \Vert \leqq  \epsilon $
for  $j = 1, \ldots , k $ .

   (ii) For any $k \in \Bbb N$, any $z \in \L_\A({X^{\otimes k}}_\A)$, any finite
subset $\{ z_1, \ldots, z_m \} \subset \L_\A({X^{\otimes k}}_\A)$ and any
$\epsilon > 0 $, there exists $q \in{\O_X}^{(0)}$ with
$\Vert q \Vert = 1$ such that
$\Vert qzq \Vert \geqq \Vert z \Vert - \epsilon $, \
$\Vert qz_i\sigma ^j(q) \Vert \leqq \epsilon $
\quad and \quad $\Vert \sigma ^j(q)z_iq \Vert  \leqq \epsilon $
for $i = 1, \ldots, m$ and $j = 1, \ldots, k $.

{\bf Remark 3.}
   We can also provide a proof of the simplicity of the $C^*$--algebra
$\Cal O _{\Lambda}$ associated with a certain class of general subshifts
studied by Matsumoto \cite{Ma} and in particular $\Cal O_{\beta}$
associated with $\beta$--shifts studied in 
\cite{KMW} by our argument.

\heading
5. Applications of the Simplicity Criterion
\endheading

\noindent{\bf  a) Real or Pseudoreal Bimodules}

 Jones Index theory for $C^*$--subalgebras
was developed  in
\cite{Wa,KW1} based on the idea that
  if $\A \subset \B$ is an inclusion of unital $C^*$--algebras with finite
index, $\B$ should be regarded
as a finite projective module over $\A$. More specifically,
   a conditional expectation $E : \B \rightarrow \A$ has finite index
if there is a finite subset $\{u_1,\ldots ,u_n \} \subset \B$, called a basis,
such that $x = \sum_{i=1}^n u_iE(u_i^*x)$ for all $x \in \B$.
   In \cite{Wa}, $\{ (u_1,u_1^*), \ldots , (u_n,u_n^*) \}$ is called a
quasi--basis.
   The index of $E$ is defined by Index $E = \sum_i u_iu_i^*$.
   It is shown that Index $E$  does not depend on the choice of basis and
Index $E$ is in the center of $\B$.
   We see that Index $E \geqq I$ and Index $E = I$ if and only if $\A =  \B$.
   We consider $X = \B_\A$ as a right Hilbert $\A$--module with the $\A$--valued
inner porduct $(x \vert y)_\A= E(x^*y)$.
 Define a unital isometric $^*$--homomorphism
$\phi : \A \rightarrow \K = \Cal K_\A(X_\A) = \Cal L_\A(X_\A)$ by $\phi (a)x = ax$
for $a \in \A$ and $x \in X = \B$.

More generally, in \cite{KW1}
the authors considered the concept of Hilbert bimodules of finite type.
We denote by $d(X)$ the (minimal) dimension of a Hilbert bimodule of finite
type.

Let us regard    $X$ as an 
object of the tensor $C^*$--category ${\H_{\A}}^b$ with objects
Hilbert $\A$--bimodules  and arrows between two objects $X$ and $Y$
the subspace $_\A\L_\A(X_\A, Y_\A)$ of 
$\L_\A(X_\A, Y_\A)$ of operators that commute with
the left $\A$--action. Assume that $X$ has finite intrinsic dimension 
in ${\H_{\A}}^b$ in the sense of \cite{LR}, i.e. there are arrows 
$R\in _{\A}\L_\A(\A_{\A}, Y_\A\otimes_{\A}X_{\A})$ and 
$\ov{R}\in _{\A}\L_\A(\A_{\A}, X_\A\otimes_{\A}Y_{\A})$ satisfying the 
conjugate equations 
$${\ov{R}}^*\otimes 1_X\circ 1_X\otimes R=1_X\ ,$$
$${R}^*\otimes 1_Y\circ 1_Y\otimes \ov{R}=1_Y\ .$$
 $X$ is called real (resp. pseudoreal) if 
there  is a solution of the conjugate equations 
 with $Y=X$ and $\ov R=R$ ($\ov R=-R$ resp.).

Let  $\A$ be a  unital $C^*$--algebra, and let $X$ be  a Hilbert bimodule over $\A$
such that $X_\A$ is finite projective.  It is then straightforward to check that
$X$ is  real  (resp. pseudoreal) if and only if
 there is an invertible $F: X\to X$ such that
$$F(axa')={a'}^*F(x)a^*\ ,$$
for $a, a'\in \A$, $x\in X$,
and satisfying  $F^2=I$ ($F^2=-I$ resp.).

The   bimodule arising from an inclusion of $C^*$--algebras $\A\subset \B$ with finite
Jones index in the sense of \cite{Wa}
 is a typical example of real Hilbert $\A$--bimodule.

\proclaim{5.1. Theorem}  Let $X$ be a Hilbert
$\A$--bimodule, with real or pseudoreal structure
defined by an element $R$ such that $\|(R^*R)^{-1}\|<1$. Then the map $\J\to\J\cap \A$
is a bijective correspondence between closed ideals of $\O_X$
and $X$--invariant ideals of $\A$. In particular $\O_X$ is simple if and
only if $\A$ is $X$--simple.
Furthermore, if $\A$ is nonnuclear, then $\Cal O_X$ is nonnuclear.
\endproclaim
\demo{Proof}
 By Theorem 4.3  we need  to show that every closed $X$--invariant ideal $J$
of $\A$ satisfies $J_X=J$ and that $X$ is $(II)$--free. 

To prove the first claim it is enough to show that for some
 $s\in{\Bbb N}$, $X^{\otimes s}$ contains
a central element element $S$ (i.e. $aS=Sa$, $a\in \A$) satisfying $(S\vert S)_\A=I$, since, then,
if $a\in J_X$ we have $a=(S\vert Sa)_\A=(S\vert \phi(a)S)_\A\in J$.

We set  $S:=R(R^*R)^{-1/2}\in X^{\otimes 2}.$ 
We have that 
 $$S^*\sigma(S)=(R^*R)^{-1/2}\sigma((R^*R)^{-1/2})\in\phi(\A)'\cap\L_\A(X_\A),$$ hence
$\|S^*\sigma(S)\|\leq\|(R^*R)^{-1}\|<1$. 
We check that the isometries 
$S_n:={\sigma}^{n-1}(S)\dots\sigma(S)S\in \A'\cap X^{2n}$ satisfy the required
properties for $(I)$--freeness. Now $(4.1)$ is a consequence
of the fact that $S_n$ is an isometry. We prove $(4.2)$.
 ${S_n}^*\sigma^n({S_n})={S_{n-1}}^*\sigma^{n-1}({S_{n-1}})\sigma^{n-1}(S^*\sigma(S))$
hence $\|{S_n}^*\sigma^n({S_n})\|<1$. 
 We now show that $X$ is $(II)$--free. 
Let  $J$ be a  closed $X$--invariant ideal of $\A$. For any $n\in{\Bbb N}$
let $\pi_n: X^{\otimes n}\to ({{X/X_J}})^{\otimes n}$ denote the n--th tensor power of the
quotient map $\pi: X\to X/X_J$. Then the elements
  $\pi_{2n}(S_{2n})\in ({X/X_J})^{\otimes 2n}$ are still isometries in 
$(\A/J)'\cap ({X/X_J})^{\otimes 2n}$ that
 satisfy the required properties.

   Finally if $\A$ is nonnuclear, then $\Cal O_X$ is nonnuclear,
because there exists a conditional expectation
    $E_\A^{\Cal O_X} : \Cal O_X \rightarrow \A$
by Lemma 2.6.
\enddemo

Note that if $d(X)$ is the minimal dimension of $X$ then $d(X)^{-1}\leq\|(R^*R)^{-1}\|$.
However, the  condition $d(X)>1$ is not enough to hold the Theorem true
in general (if we do not assume that $\A$ has trivial center)
 as the following counterexample shows. 

{\bf Example 1.}
 Let $\A_1 = \Bbb C = \B_1$ and $\A_2 = \Bbb C \subset B_2 = \Bbb C^2$.
Consider conditional expectations $E_1 =id :\B_1 \rightarrow \A_1$ and 
$E_2 : \B_2 \rightarrow \A_2$ defined by 
$E_2(x,y) = (\frac{x+y}{2},\frac{x+y}{2})$.
   Set $\A = \A_1 \oplus \A_2 \subset \B = \B_1 \oplus \B_2$ and 
$E=E_1\oplus E_2: \B\rightarrow \A$
$X_1 = (\B_1)_{\A_1}$ and 
$X_2 = (\B_2)_{\A_2}$,
 $X = \B_\A$ with  inner products defined by 
$E$,
and left and right actions of $\A$ defined by multiplication. 
Since $X = X_1 \oplus X_2$, we have 
$$
   \Cal O_X \cong \Cal O_{X_1} \oplus \Cal O_{X_2} 
\cong {\Bbb C}({\Bbb T}) \oplus \Cal O_2
$$
Thus it is clear that $\O_X$ is not simple. We have 
$d(X)>1$  but there is no  bijective  correspondence between ideals, because $\O_1$ is
isomorphic to ${\Bbb C}({\Bbb T})$. In fact $X$ is not $(II)$--free since $X/X_{\A_2}={\Bbb C}$
is not $(I)$--free.

However, if we assume that $\A$ is $X$--simple, the conclusions of the above Theorem hold
 in the following important particular case.

\proclaim{5.2. Theorem}
   Let $\A \subset \B$ be an inclusion of unital $C^*$--algebras.
Let $E : \B \rightarrow \A$ be a conditional expectation with
finite index in the sense of \cite{Wa}.
   Consider $X = \B_\A$ as a right Hilbert $\A$--module with
$(x \vert y)_\A = E(x^*y)$ for $x,y \in \B$ and define 
   a left $\A$--action $\phi : \A \rightarrow \Cal K_\A(X_\A) $    by
$\phi(a)x = ax \quad $  for $a \in \A$ and $x \in X$.
   If $\A$ is $X$--simple and Index $E \not= I$, then the
$C^*$--algebra $\Cal O_X$ is simple.
\endproclaim
\demo{Proof}
  We show
that $X$ is $(I)$--free, so the result is a consequence of Corollary 4.4.
 Since Index $E \not= I$, we have that the Jones projection $e_\A \not= I$ and there exists
an element $x_0 \in X$ such that $e_\A(x_0) \ne x_0$
   Since $e_\A \in \A' \cap \Cal K_\A(X_\A)$, we can define
$q_k = e_\A \otimes \cdots \otimes e_\A \otimes (1-e_\A) \in
\Cal K_\A({X^{\otimes k+1}}_\A)$ for each $k \in \Bbb N$.
   Then it is easy to see that $q_k \in \A'\cap\O_X$ and for any $m \in \Bbb N$
with $1 \leqq m \leqq k$, we have $q_k\sigma ^m(q_k) = 0$.
   Since
$$
\align
 & \quad \quad (q_k(1\otimes \cdots \otimes 1 \otimes x_0)\ \vert \
                    1\otimes \cdots \otimes 1 \otimes (x_0-e_\A(x_0)))_\A \\
 & \quad =         (1\otimes \cdots \otimes 1 \otimes (x_0-e_\A(x_0)) \
            \vert \ 1\otimes \cdots \otimes 1 \otimes (x_0-e_\A(x_0)))_\A \\
 & \quad = (x_0 - e_\A(x_0) \ \vert \ x_0 - e_\A(x_0))_\A \ne 0
\endalign$$
   Therefore $q_k \ne 0$.
   Let $J = \{ a \in \A ; ab = aE(b) \quad \text{for all} \quad b \in \B\}$.
Then $J$ is a closed ideal of $\A$ and it is $X$--invariant.  In fact for
$a \in J$, $\quad x, y, b \in \B$ we have
$$
   (x\vert ay)_\A b = E(x^*ay)b = E(x^*aE(y))b 
= E(x^*)aE(E(y)b) =
 (x\vert ay)_\A E(b).
$$
Thus $J=\{0\}$ since $I\notin J$ and $\A$ is $X$--simple.
   Now the
$^*$--homomorphism $\A \owns a \rightarrow \phi(a)q_k \in 
{\O_X}^{(0)}$ has kernel exactly $J$,
hence it is faithful. This shows that
    $X$ is $(I)$--free.

\enddemo

{\bf Example 2.}
\quad Let $\A$ be a unital $C^*$--algebra, $G$ a finite group and
$\alpha : G \rightarrow \Aut\A$ an automorphic action. 
Consider the crossed product $\B = \A \rtimes _{\alpha} G$
and the  usual conditional expectation $E : \B \rightarrow \A$ and let
$\lambda _g \in \B$ be the unitary in $B$ implementing $\alpha_g$. 
   For the inclusion $\A \subset \B$, let $X = \B_\A$ be a right Hilbert
module with $\phi : \A \rightarrow \Cal K_\A(X_\A)$ as above.
If the order $n$ of G is not equal to one and $\A$ is simple,
then $\Cal O_X$ is simple.
   The $C^*$--algebra $\Cal O_X$ is uniquely determined by the
generators $\A$ and
n isometries $\{ S_g ; g \in G \} $ with the relations
$aS_g = S_g\alpha _g^{-1}(a)$ and $ \sum _g S_gS_g^* = I$.
   A concrete realization is given in $\B \otimes \Cal O_n$.
Let $\{ T_g ; g \in G \}$ be  the generators of $\Cal O_n$. 
Then $\Cal O_X$ is generated by $\{ a\otimes I ; a \in \A \}$ and
$\{ S_g = \lambda _g \otimes T_g \}$.
   If the action $\alpha$ is trivial, then $\Cal O_X$ is isomorphic
to the tensor product $\A \otimes \Cal O_n$. 
   Therefore  $\Cal O_X$ is looks like a twisted tensor
product of $\A$ by the Cuntz algebra $O_n$, that is ,
"$\Cal O_X \cong \A\otimes_{\alpha} \Cal O_n$"  by an abuse of notation.
  
{\bf Example 3.}
\quad Let $G$ be a free group of two generators and $H$ a subgroup
of $G$ with $[G:H] = n$.  Then $H$ is also a free group of $n+1$ generators.
By  \cite{Po} we  have an inclusion of simple $C^*$--algebras
$\A = C^*_r(H) \subset \B = C^*_r(G)$ and a conditional expectation
$E : \B \rightarrow \A$ such that
   $E(\sum_{g\in G}x_g\lambda _g) = \sum_{h\in H}x_h\lambda _h$.  
Put $X = \B_\A$ as above.  Then $\Cal O_X$ is simple and nonnuclear
by Theorem 5.1.

{\bf Example 4.}
  \quad  Let $e_1,e_2,\ldots $  be a sequence of Jones projections such that
$e_ie_{i\pm 1}e_i = (4\cos^2\pi /5)^{-1}  e_i$ for $i = 1,2,\ldots $.
   Let $\B = C^*\{ 1,e_1,e_2,e_3,\ldots \} \supset
\A = C^*\{ 1,e_2,e_3,\ldots \} $ and $E : \B \rightarrow \A $ a conditional
expectation of finite index as in  [Wa; Example 3.3.14]. 
   Since $\A$ is simple and Ind $E = 4\cos^2\pi /5$, we have that
 $\Cal O_X$ is simple by Theorem 5.1.

\noindent{\bf b) Cuntz--Krieger Bimodules}

In this subsection we focus our attention on a class of finite projective
Hilbert bimodules over finite direct sum of simple unital
$C^*$--algebras. We
start with a couple of easy results that we shall need.\medskip

\proclaim{5.3. Proposition}  If $X$ is a full Hilbert $C^*$--bimodule over
a unital $C^*$--algebra $\A$
then there is a multiplet $\hat y=\{y_1,\dots, y_n\}$ in $X$ such that
$$\sum_i (y_i\vert y_i)_\A=I\ .$$
The c. p. map $\Phi_{\hat y}: T\in\O_X\to\sum_i {y_i}^*T{y_i}
\in\O_X$ satisfies
$$\Phi_{\hat y}(A\sigma(B))=\Phi_{\hat y}(A)B\ ,
 A\in\O_X\ , B\in{\A}'\cap\O_X\ .$$
In particular, $\Phi_{\hat y}$ is a left inverse of $\sigma$ on
$\A'\cap\O_X$.
\endproclaim
\demo{Proof}
Since $X$ is full and $\A$ is unital, there are elements $z_i$ and $z'_i$
such that $\|\sum_i (z_i\vert z'_i)_\A-I\|<1$. In particular
 $\sum_i (z_i\vert z'_i)_\A$ is invertible,
so we may assume $\sum_i (z_i\vert z'_i)_\A=1$. The matrix
$(\theta_{z'_i, z'_j})$ is positive, thus there are $T_{i, k}\in\L_\A(X_\A)$
 such that
$(\theta_{z'_i, z'_j})=\sum_k T_{i, k}{T_{j, k}}^*$. It follows that
$$\sum_{i, j, k}({T_{k, i}}^*z_i\vert {T_{j, k}}^*z_j)_\A =I\ .$$
The elements  $y_k=\sum_j {T_{j, k}}^*z_j$ satisfy the desired properties.
\enddemo

\proclaim{5.4. Lemma} If $\A$ is unital and $X$ is a full Hilbert
bimodule over $\A$
then $$\Z(\L_\A(X_\A))=\{x\in X\to xa\in X\ ,\quad a\in\Z(\A)\ \}$$
\endproclaim
\demo{Proof} Clearly the set at the right hand side is contained in the center
of $\L_\A(X_\A)$.
 Let
 $y_1\ ,\dots , y_n$ be a multiplet in $X$ such that
$\sum_i (y_i\vert y_i)_\A=I$.
 If $T\in\Z(\L_\A(X_\A))$  then for any $x\in X$,
$$Tx=x\sum_i(y_i\vert Ty_i)_\A\ .$$ Replacing $x$ by $xu$, with $u$ unitary in
$\A$,
shows  that $a=\sum_i(y_i\vert Ty_i)_\A$ is in the center of $\A$, 
and the proof is complete. 
\enddemo

 We denote by $r_n: a\in\Z(\A)\to r_n(a)\in\Z(\L_\A({X^{\otimes n}}_\A))$ the action of
$\Z(\A)$ on $X^{\otimes n}$ given by right multiplication. We write $r$ for $r_1$.

Let $\A$ be the direct sum of $d>1$ unital simple
$C^*$--algebras,
$$\A=\A_1\oplus\dots\oplus \A_d$$
and let $p_1\ ,\dots, p_d$ denote the minimal central projections of
$\A$.
We associate to $X$ the  matrix $A\in M_d(\{0, 1\})$ by
$$A(i, j)=0\ ,\quad  \hbox{\rm if}\  \phi(p_i)r(p_j)=0$$,
$$A(i, j)=1\ ,\quad  \hbox{\rm if}\  \phi(p_i)r(p_j)\ne 0\ .$$
Recall that $A$ is called {\it irreducible} if for any pair of $i, j$
there is $n\in{\Bbb N}$ such that $A^n(i, j)>0$.\medskip

\proclaim{5.5. Proposition} $\A$ is $X$--simple if and
only if the associated matrix $A$ is irreducible.
\endproclaim
\demo{Proof} $\A$ is $X$--simple if and only if for all $j$,
the smallest $X$--invariant ideal containing $p_j$ is $\A$,
that is to say, for all $i$ there is $n$ such that
$$\{(x\vert \phi (p_j)yp_i)\ ,\quad x, y\in X^{\otimes n}\}\ne 0.$$
In other words, $\phi (p_j)r_n(p_i)\ne 0$.
To complete the proof we need to show that
$\phi (p_i)r(p_{i_1})\ne 0$, $\dots$,  
$\phi (p_{i_n})r(p_j)\ne 0$ imply
$\phi (p_i)r(p_{i_1})\dots r_n(p_j)\ne 0$ and  hence
$\phi (p_i)r_n(p_j)\neq 0$. By induction,
assume that   
$\phi (p_i)r(p_{i_1})\dots r_{n-1}(p_{i_{n}})\ne 0$.
The map
$$T\in\L_\A({X^{\otimes n}}_\A)\to r_{n}(p_j)T$$
is injective on $r_{n-1}(p_{i_n})\L_\A({X^{\otimes n}}_\A)$
so the proof is completed by the previous Lemma.
\enddemo

Following Cuntz and Krieger, we consider the commutative  $C^*$--subalgebra
$\D_{A}$ generated
by $r_k(p_i)$, $i=1,\dots, d$, $k=0, 1, 2,\dots$.
 \medskip

\proclaim{5.6. Proposition} The spectrum of $\D_A$ is the compact
space
$$\xi_{A}=\{(x_k)_{k\in{\Bbb N}}, x_k=1,\dots, d, A(x_{k}, x_{k+1})=1\}\ .$$
Each projection $p_i$ corresponds to the characteristic function of
$$Z(i)=\{(x_k)\in \xi_A, x_1=i\}.$$
Furthermore
the endomorphism $\sigma$ acts on $\xi_A$ as the one--sided shift on the left.
\endproclaim
\demo{Proof} The spectrum of $\D_A$ is given by
$$\{(x_k)_{k\in{\Bbb N}}, x_k=1,\dots, d, p_{x_1}r(p_{x_2})\dots
r_n(p_{x_{n+1}})\ne 0, n\in{\Bbb N}\}\ ,$$
 so the arguments of the
previous Proposition complete the proof. 
\enddemo

\proclaim{5.7. Theorem}  Let $\A$ be the direct sum of $d>1$ unital simple 
$C^*$--agebras as above.
Let $X$ be a finite projective Hilbert
bimodule over $\A$. If the associated matrix $A$ satisfies condition
$(I)$ of \cite{CK} then $X$ is $(I)$--free. In particular, if $A$ is irreducible
and non permutation
then  $\O_X$ is a simple
$C^*$--algebra.
\endproclaim
\demo{Proof} We can apply the arguments of [CK; Lemmas 2.6, 2.7] 
to  $\L_\A({X^{\otimes k}}_\A)$ in place of $\F_k$, and to a basis 
$\{x_i\}$ of $X$ in
place of the partial isometries $S_i$
and deduce that for all $k\ , r\in{\Bbb N}$ there is a projection
$q\in \D_{A}$  such that
$$\|{\sigma}^r(q)T\|=\|T\|\ , \quad T\in \L_\A({X^{\otimes r}}_\A)\ ,$$
$$q{\sigma}^j(q)=0\ ,\quad 0<|j|<k\ .$$
Approximating $q$ by a projection  in some $\Z(\L_\A({X^{\otimes p}}_\A))$
gives an element with the desired properties.
\enddemo

\heading
6. Bimodules over Purely Infinite $C^*$--Algebras
\endheading

 In this section we present a class of examples of purely infinite 
simple $C^*$--algebras  of the form $\Cal O _X$.

\proclaim{6.1. Theorem}  Let $X$ be a finite projective right Hilbert
bimodule over a unital  purely infinite simple $C^*$--algebra $\A$. 
If $X$ is $(I)$--free, then
$\O_X$ is simple and purely infinite.
\endproclaim
\demo{Proof}
Note that if $Y$ is a finite projective Hilbert module over $\A$ then
 $\L_\A({Y}_\A)$ is isomorphic to a full
corner of some $M_d(\A)$, so it is simple and purely infinite. Now
${\O_X}^{(0)}$ is  simple and purely infinite 
since it is the inductive limit of $\L_\A({X^{\otimes r}}_\A)$.
   
   Let $T$ be a non zero positive element in $\Cal O_X$.  Since $m_0(T)$
is a non zero positive element in the purely infinite $C^*$--algebra 
${\O_X}^{(0)}$, there exists $W \in {\O_X}^{(0)}$ such that
$W^*m_0(T)W = 1$.  Put $S = W^*TW$. Then $m_0(S) = 1$.  
   Choose $0 < \varepsilon < 1/5$.
   Then there exists $B \in {}^0\Cal O_X$ with $\|S - B \| < \varepsilon$.
   Put $B_0 = m_0(B) \in {^0\O_X}^{(0)}$.
Then $\| 1 - B_0 \| = \| m_0(S - B) \| < \varepsilon$.
  Let $C = B - B_0 + I \in {}^0\Cal O_X$.  
  Then $\|C - B \| < \varepsilon $ and $\|S - C \| < 2\varepsilon $.
  Consider the expansion  $C = \sum _{j=-k}^k C_j \in {}^0\Cal O_X$,
with $C_j \in \L_\A({X^{\otimes p}}_\A, {X^{\otimes p+j}}_\A)$ for some $p$.
  Note that $C_0 = m_0(C) = I$.
  By the proof of Theorem 4.3, considering the product of operators
that appear in the iterating argument, $(I)$--freeness provides an
element $R \in \A'\cap {}^0\Cal O_X$ satisfying
$R^*R \in {\O_X}^{(0)}$, $\|R^*CR - R^*C_0R \| < \varepsilon$
and $\|R^*R\| = \|R^*C_0R\| = \|C_0\| = 1 $.
  Then
$$
  \|R^*SR - R^*R \| \le \|R^*(S - C)R\| + \|R^*(C - C_0)R\| < 3\varepsilon
$$
  Since $R^*R$ has norm one,  there exists $V \in {\O_X}^{(0)}$
such that $V^*R^*RV = I$ and $\|V\| < 1 + \varepsilon$.
   Then we have 
$$
 \|V^*R^*SRV - I \| = \|V^*R^*SRV - V^*R^*RV \| 
  < 3\varepsilon (1+\varepsilon  )^2 < 108/125 < 1
$$
  Hence $V^*R^*SRV$ is positive invertible and there exists $Y \in \Cal O_X$
such that $Y^*TY = 1$.  Therefore $\Cal O_X$ is purely infinite. 
\enddemo               

  {\bf Remark 1.} In the above Theorem, instead of assuming $\A$ is purely
infinite, we may assume that $\A$ has real rank zero and for any natural 
number 
$r$ and for any non zero projection $P \in \Cal L_\A(X^{\otimes r}_\A)$ there
exists $x \in X^{\otimes r}$ with $Px = x$ and $(x\vert x)_\A = 1$.
Just note the fact that $x^*x = I$ and  
$$
   I \geq x^*Px \geq x^*\theta_{x,x}x = x^*xx^*x = I
$$ 

{\bf Remark 2.} If we know a priori that $\Cal O_X$ is simple  
 then we can also conclude, without  assuming that $X$ is $(I)$--free, that
$\Cal O_X$ is purely infinite as follows.
By [Pe; Theorems 8.10.4; 8.10.10] ${\O_X}^{(0)}$ has the trivial 
relative  commutant in $\O_X$.
Let now $x$ be the image in $\O_X$ of a non zero element in $X$.
By the pure infiniteness of
${\A}$ there is $a\in{\A}$ such that $a^*x^*xa=I$, so replacing $x$ by
$xa$, we may assume that $x$ is an isometry. Note that $\O_X$ is
generated by ${\O_X}^{(0)}$ and $x$.
Now the triviality of the relative commutant implies that all the powers
of the endomorphism $\rho$ implemented by $x$ on ${\O_X}^{0}$ are 
outer. If $x$ is unitary then $\O_X$
is the crossed product of ${\O_X}^{0}$ by the  automorphism $\rho$
  so it is simple \cite{Ki} and purely infinite
 [JKO; Theorem 3.2]. If $x$ 
is a proper
isometry, then $\O_X$ is the crossed product of a simple purely infinite
$C^*$--algebra by a proper corner endomorphism, such that all powers are
outer
so by
[JKO; Theorem 2.1], $\O_X$ is simple and purely infinite.

\noindent{\bf Acknowledgment} 
Part of this paper was written during a visit of C.P. to 
the Matematisk Institut of Kobenhavn. She
 wishes to thank  Ryszard Nest, Hiroyuki Osaka, Gert K. Pedersen, Carl Winsl\o w for fruitful discussions
and for warm hospitality.
\heading
References
\endheading

\item{[AEE]}B.~Abadie, S.~Eilers and R.~Exel,
{\it Morita equivalence for crossed products by Hilbert $C^*$--bimodules\/},
Preprint (1995).

\item{[BMS]} L.~G.~Brown, J.~Mingo and N.~Shen, 
{\it Quasi--multipliers and embeddings
of Hilbert $C^*$--bimodules\/},  Canad.~J.~Math.,  {\bf 46, No 6} (1994), 
1150--1174.

\item{[Co1]} A.~Connes: {\it Une classification des facteurs de type III\/},
Ann.~Scient.~\`Ec. Norm.~Sup., {\bf 6} (1973), 133--252.

\item{[Co2]} A.~Connes: {\it Outer conjugacy classes of automorphisms
of factors\/}, Ann.~Scient.~\`Ec. Norm.~Sup., {\bf 8} (1975), 383--420.

\item{[CT]} A.~Connes and M.~Takesaki: {\it The flow of weights on factors
of type III}, Tohoku~Math.~J., {\bf 29} (1977), 473--575.

\item{[Cu]}J.~Cuntz,
{\it Simple $C^*$--algebras generated by isometries\/},
Comm.~Math.~Phys., {\bf 57} (1977), 173--185.

\item{[CK]}J.~Cuntz and W.~Krieger,
{\it A class of $C^*$--algebras and topological Markov chains\/},
Invent.~Math., {\bf 56} (1980), 251--268.

\item{[Cu2]}J.~Cuntz, 
{\it A class of $C^*$--algebras and topological Markov chains II:
 Reducible chains and the Ext--functor for $C^*$--algebras\/},
Invent.~Math., {\bf 63} (1981), 25--40.

\item{[DR1]}S.~Doplicher and J.~E.~Roberts, 
{\it Duals of compact Lie groups
realized in the Cuntz algebras and their actions on $C^*$--algebras\/},
J.~Funct.~Anal., {\bf 74} (1987), 96--120.

\item{[DR2]}S.~Doplicher and J.~E.~Roberts, 
{\it Endomorphisms of $C^*$--algebras,
crossed products and duality for compact groups\/}, Ann.~Math., {\bf 130} (1989),
 75--119.

\item{[DR3]}S.~Doplicher and J.~E.~Roberts, 
{\it A new duality theory for compact groups\/},
Invent.~Math., {\bf 98} (1989), 157--218.

\item{[DPZ]}S.~Doplicher, C.~Pinzari and R.~Zuccante, 
{\it On the $C^*$--algebra
generated by a Hilbert bimodule\/}, preprint.

\item{[El]}G.~Elliott, 
{\it  Some simple $C^*$--algebras constructed as crossed products with
discrete outer automorphism groups\/},
Publ.~RIMS, Kyoto Univ., {\bf 16} (1980), 299--311.

\item{[Ex]}R.~Exel,
{\it Circle actions on $C^*$--algebras, partial automorphisms and
a generalized Pimsner--Voiculescu exact sequence\/},
J.~Funct.~Anal., {\bf 122} (1994), 361--401.

\item{[GHJ]}F.~M.~Goodman, P.~de~la~Harpe and V.~F.~R.~Jones,
{\it Coxeter graphs
and towers of algebras\/}, Springer--Verlag, New York, (1989).

\item{[I1]}M.~Izumi, 
{\it Subalgebras of infinite $C^*$--algebras with
finite Watatani indices. I. Cuntz algebras\/},
Comm.~Math.~Phys., {\bf 155} (1993), 157--182.

\item{[I2]}M.~Izumi, 
{\it Subalgebras of infinite $C^*$--algebras with
finite Watatani indices. II. Cuntz--Krieger algebras\/}, 
Preprint.

\item{[I3]} M.~Izumi, 
 {\it Index Theory of Simple $C^*$--algebras\/}, Talk
given at the Fields Institute,  1995.

\item{[JKO]}J.~A~Jeong, K.~Kodaka and H.~Osaka, 
{\it Purely infinite simple 
$C^*$--crossed products I\/},  Canad.~Math.~Bull.~Vol., to appear.

\item{[Jo]}V.~F.~R.~Jones, 
{\it Index for subfactors\/},
Invent.~Math., {\bf 72} (1983), 1--25.

\item{[KW1]}T.~Kajiwara and Y.~Watatani, 
{\it Jones index theory by Hilbert $C^*$--bimodules and K--theory\/},
preprint (1996).

\item{[KW2]}T.~Kajiwara and Y.~Watatani, 
{\it Crossed products of Hilbert $C^*$--bimodules by countable
discrete groups\/},
to appear in Proc.~Amer.~Math.~Soc. 

\item{[KWP]}T.~Kajiwara, Y.~Watatani and C.~Pinzari, {\it Hilbert $C^*$--bimodules
and countably generated Cuntz--Krieger algebras\/}, in preparation. 
                                          
\item{[Ka]}Y.~Katayama, 
{\it Generalized Cuntz algebras $\Cal O_N^M$\/},
RIMS~Ko\-kyu\-ro\-ku, {\bf 858}\linebreak (1994), 131--151.

\item{[KMW]}Y.~Katayama, K.~Matsumoto and Y.~Watatani, 
{\it Simple $C^*$ algebras arising from $\beta$--expansion of
real numbers\/}, Preprint (1996)

\item{[Ki]}A.~Kishimoto,
{\it Outer automorphisms and reduced crossed products of simple
$C^*$--algebras\/},
Comm.~Math.~Phys., {\bf 81} (1981), 429--435.

\item{[LR]}R.~Longo and J.~E.~Roberts, 
 {\it A theory of dimension\/}, to appear in K--theory.

\item{[Ma]}K.~Matsumoto,
{\it On $C^*$--algebras associated with subshifts\/},
to appear in Internat.~J.~Math. 

\item{[MS]} P.~S.Muhly and B.Solel, {\it On the simplicity of some
Cuntz--Pimsner algebras\/}, preprint (1996).

\item{[Ol]} D.~Olesen: {\it Inner automorphisms of simple $C^*$--algebras\/},
Comm.~Math.~Phys., {\bf 44} (1975), 175--190.

\item{[OP1]} D.~Olesen and G.~K.~Pedersen: {\it Applications of the Connes
spectrum to $C^*$--dynamical systems\/}, 
 J.~Funct.~Anal., {\bf 30} (1978), 179--197.

\item{[OP2]} D.~Olesen and G.~K.~Pedersen: {\it Applications of the Connes
spectrum to $C^*$--dynamical systems III\/}, 
 J.~Funct.~Anal., 
{\bf 45} (1982), 357--390 .

\item{[Pe]}G.~K.~Pedersen, 
{\it $C^*$--al\-ge\-bras and their auto\-mor\-phism groups\/},
Academic \linebreak Press, New York (1990).

\item{[PP]}M.~Pimsner and S.~Popa, 
{\it Entropy and index for
subfactors\/},  Ann.~Sci.~Ec. Norm. Sup., {\bf 19} (1985),  57--106.

\item{[Pi]}M.~Pimsner, 
{\it A class of $C^*$--algebras generalizing both Cuntz--Krieger algebras and
Crossed products by $\Bbb Z$\/},
preprint.

\item{[P]}C.~Pinzari, {\it The ideal structure of Cu\-ntz--Krie\-ger\---Pimsner algebras
and Cun\-tz\---Krieger algebras over infinite matrices\/}, to appear in 
the proceedings of the Rome conference ``Operator algebras and Quantum Field Theory'',
July 1996.

\item{[Po]}R.~Powers,
{\it Simplicity of $C^*$--algebra associated with the free group on two
generators\/},
Duke~Math.~J., {\bf 42} (1975), 151--156.

\item{[Ri]}M.~Rieffel, 
{\it Induced representations of $C^*$--algebras\/},
Adv.~Math., {\bf 13} (1974), 176--257.

\item{[Ro]}J.~E.~Roberts, 
{\it Cross products of von Neumann algebras by group duals\/},
Symp.~Math., {\bf 20} (1976), 335--363.

\item{[R\o1]}M.~R\o rdam, 
{\it Classification of Cuntz--Krieger algebras\/},
K--theory {\bf 9} (1995), 31--58.

\item{[R\o2]}M.~R\o rdam, 
{\it Classification of certain infinite simple $C^*$--al\-ge\-bras\/},
J.~Fun\-ct.~Anal., {\bf 131} (1995), 415--458.

\item{[Wa]}Y.~Watatani, 
{\it Index for $C^*$--sub\-al\-ge\-bras\/},
Memoir~Amer.~Math.~Soc., {\bf 424} (1990).

\end{document}